\documentclass{amsart}

\usepackage{amssymb}
\usepackage[all]{xy}

\def\Z{{\mathbb Z}}
\def\Q{{\mathbb Q}}
\def\C{{\mathbb C}}
\def\N{{\mathbb N}}
\def\P{{\mathbb P}}
\def\F{{\mathbb F}}

\def\A{{\mathcal A}}
\def\B{{\mathcal B}}
\def\cC{{\mathcal C}}
\def\cG{{\mathcal G}}
\def\J{{\mathcal J}}
\def\L{{\mathcal L}}
\def\M{{\mathcal M}}
\def\O{{\mathcal O}}
\def\cR{{\mathcal R}}

\def\U{{\mathcal U}}

\def\rhotilde{{\tilde{\rho}}}
\def\rhobar{{\overline{\rho}}}
\def\rhohat{{\hat{\rho}}}
\def\nutilde{{\tilde{\nu}}}

\def\G{\Gamma}
\def\Ghat{\widehat{\G}}
\def\Gtilde{\widetilde{\G}}
\def\Mtilde{\widetilde{M}}

\def\Mbar{\overline{\M}}
\def\Qbar{\overline{\Q}}

\def\dot{\bullet}

\def\l{\ell}
\def\Ql{{\Q_\ell}}
\def\Zl{{\Z_\ell}}

\def\Fl{{\F_\ell}}

\def\Cbar{\overline{C}}
\def\Tbar{\overline{T}}

\def\phibar{\overline{\varphi}}

\def\Gm{{\mathbb{G}_m}}
\def\Ga{{\mathbb{G}_a}}
\def\Sp{{\mathrm{Sp}}}
\def\SL{{\mathrm{SL}}}
\def\GL{{\mathrm{GL}}}
\def\GSp{{\mathrm{GSp}}}

\def\prol{{(\ell)}}
\def\un{\mathrm{un}}
\def\arith{\mathrm{arith}}

\def\cts{\mathrm{cts}}
\def\ur{\mathrm{ur}}

\def\x{{\mathbf{x}}}
\def\vecs{\vec{\mathbf{v}}}

\def\Gp{\mathsf{Group}}
\newcommand\Vect{\operatorname{\mathsf{Vec}}}

\newcommand\im{\operatorname{im}}

\newcommand\Hom{\operatorname{Hom}}

\newcommand\Spec{\operatorname{Spec}}
\newcommand\Diff{\operatorname{Diff}}
\newcommand\Aut{\operatorname{Aut}}
\newcommand\Out{\operatorname{Out}}
\newcommand\Gr{\operatorname{Gr}}
\newcommand\Heis{\operatorname{Heis}}
\newcommand\ord{\operatorname{ord}}
\newcommand\Jac{\operatorname{Jac}}
\newcommand\supp{\operatorname{supp}}
\newcommand\Ind{\operatorname{Ind}}

\newtheorem{theorem}{Theorem}[section]
\newtheorem{lemma}[theorem]{Lemma}
\newtheorem{proposition}[theorem]{Proposition}
\newtheorem{corollary}[theorem]{Corollary}
\newtheorem{bigtheorem}{Theorem}
\newtheorem{bigcorollary}[bigtheorem]{Corollary}

\theoremstyle{definition}

\newtheorem{example}[theorem]{Example}

\theoremstyle{remark}
\newtheorem{remark}[theorem]{Remark}

\begin{document}

\title{Relative Pro-$\ell$ Completions of Mapping Class Groups}

\dedicatory{To Professor Gus Lehrer on his 60th birthday}

\author{Richard Hain}
\address{Department of Mathematics\\ Duke University\\
Durham, NC 27708-0320}
\email{hain@math.duke.edu}

%%\thanks{}

\author{Makoto Matsumoto}
\address{Department of Mathematics\\ Graduate School of Science\\
Hiroshima University, 739-8526 Japan}
\email{m-mat@math.sci.hiroshima-u.ac.jp}

\thanks{Supported in part by the Scientific Grants-in-Aid 16204002 and 19204002
and by the Core-to-Core grant 18005 from the Japan Society for the Promotion of
Science; and grants DMS-0405440 and DMS-0706955 from the National Science
Foundation.}

\date{\today}

\begin{abstract}

Fix a prime number $\ell$. In this paper we develop the theory of relative
pro-$\l$ completion of discrete and profinite groups --- a natural
generalization of the classical notion of pro-$\ell$ completion --- and show
that the pro-$\l$ completion of the Torelli group does not inject into the
relative pro-$\l$ completion of the corresponding mapping class group when the
genus is at least $2$. (See Theorem~\ref{thm:mcg_ladic} below.) As an
application, we prove that when $g\ge 2$, the action of the pro-$\ell$
completion of the Torelli group $T_{g,1}$ on the pro-$\ell$ fundamental group of
a pointed genus $g$ surface is not faithful.

The choice of a first-order deformation of a maximally degenerate stable curve
of genus $g$ determines an action of the absolute Galois group $G_\Q$ on the
relative pro-$\ell$ completion of the corresponding mapping class group. We
prove that for all $g$ all such representations are unramified at all primes
$\neq \ell$ when the first order deformation is suitably chosen. This proof was
communicated to us by Mochizuki and Tamagawa.

\end{abstract}

\maketitle

\section{Introduction}

Suppose that $\G$ is a discrete or profinite group and that $P$ is a profinite
group. Suppose that $\rho : \G \to P$ is a continuous, dense
homomorphism.\footnote{A homomorphism $\phi : \G \to G$ from a group $\G$ to a
topological group $G$ is {\em dense} if its image is dense in $G$.} The {\em
relative pro-$\l$ completion $\G^{\prol,\rho}$ of $\G$ with respect to $\rho:\G
\to P$} is defined by
$$
\G^{\prol,\rho}:= \lim_{\longleftarrow} G_\phi
$$
where the limit is taken over the inverse system of commutative triangles
$$
\xymatrix{
\G \ar[rr]^\phi \ar[dr]_\rho & & G_\phi \ar[dl]^{\rho_\phi} \cr
 & P
}
$$
where $G_\phi$ is a profinite group, $\phi$ is a continuous dense homomorphism,
$\rho_\phi$ is continuous, and $\ker \rho_\phi$ is a pro-$\ell$ group. There are
natural homomorphisms $\G \to \G^{\prol,\rho} \to P$ whose composition is
$\rho$.

Relative pro-$\ell$ completion with respect to the trivial representation $\rho$
of a discrete or profinite group is simply the classical pro-$\l$ completion of
the group. Relative pro-$\ell$ completion with respect to non-trivial
representations is natural in arithmetic geometry.

Since finite $\ell$-groups are nilpotent, the classical pro-$\ell$ completion of
a group $\Gamma$ with vanishing $H_1(\G,\Fl)$ is trivial. Since all mapping
class groups in genus $g$ are perfect when $g\ge 3$, their pro-$\ell$
completions are trivial.\footnote{The pro-$\ell$ completions of mapping class
groups vanish in genus 1 when $\ell\neq 2,3$ and in genus 2 when $\ell \neq 2,5$
for similar reasons.} In contrast, the natural relative pro-$\ell$ completions
of mapping class groups are large and more closely reflect their structure. For
example, when the number of marked points is at least 1, a mapping class group
injects into its relative pro-$\ell$ completion.

Suppose that $g$, $n$, and $r$ are non-negative integers satisfying
$2g-2+n+r>0$. Denote by $\G_{g,n,\vec{r}}$ the mapping class group of a compact
oriented surface of genus $g$ with $n$ distinct points and $r$ distinct non-zero
tangent vectors (or, alternatively, $r$ boundary components) whose  anchor
points are distinct from the $n$ marked points.\footnote{A detailed definition
is given in Section~\ref{sec:mcgs}.}\footnote{We adopt the standard convention
that $n$ and $\vec{r}$ are omitted when they are zero. For example, $\G_g =
\G_{g,0,\vec{0}}$, $\G_{g,\vec{r}} = \G_{g,0,\vec{r}}$ and $\G_{g,n} =
\G_{g,n,\vec{0}}$.} Denote by $T_{g,n,\vec{r}}$ its Torelli subgroup. Denote the
group of symplectic $2g\times 2g$ matrices with entries in a ring $R$ by
$\Sp_g(R)$. Denote the relative pro-$\ell$ completion of $\G_{g,n,\vec{r}}$ with
respect to the natural homomorphism $\G_{g,n,\vec{r}} \to \Sp_g(\Zl)$ by
$\G_{g,n,\vec{r}}^\prol$. Denote the pro-$\l$ completion of $T_{g,n,\vec{r}}$ by
$T_{g,n,\vec{r}}^\prol$.

For $m\in \N$, the level $m$ subgroup $\Sp_g(\Z)[m]$ of $\Sp_g(\Z)$ is defined
to be the kernel of $\Sp_g(\Z) \to \Sp_g(\Z/m\Z)$. The level $m$ subgroup of
$\G_{g,n,\vec{r}}$ is defined to be the inverse image of $\Sp_g(\Z)[m]$ in
$\G_{g,n,\vec{r}}$. Denote the relative pro-$\ell$ completion of
$\G_{g,n,\vec{r}}[m]$ with respect to the homomorphism to the closure of its
image in $\Sp_g(\Zl)$ by $\G_{g,n,\vec{r}}[m]^\prol$.

Define a discrete group $T$ to be {\em $\l$-regular} if every finite quotient of
$T$ that is an $\l$-group is the quotient of a torsion free nilpotent quotient
of $T$. Free groups, surface groups and pure braid groups are all $\l$-regular.
We suspect, but were unable to prove, that when $\l$ is odd and $n+r>0$,
$T_{g,n,\vec{r}}$ is $\l$-regular.

\begin{bigtheorem}
\label{thm:mcg_ladic}
Suppose that $g\ge 3$ and that $r$ and $n$ are non-negative integers. For each
prime number $\ell$ and each positive integer $m$, the kernel of the natural
homomorphism $T_{g,n,\vec{r}}^\prol \to\G_{g,n,\vec{r}}[m]^\prol$ is central in
$T_{g,n,\vec{r}}^\prol$ and contains a copy of $\Zl$. If $T_{g,n,\vec{r}}$ is
$\l$-regular, then there is an exact sequence
$$
0 \to \Zl \to T_{g,n,\vec{r}}^\prol \to \G_{g,n,\vec{r}}[m]^\prol
\to \Sp_g(\Zl)[m] \to 1.
$$
\end{bigtheorem}

The theorem is trivially false in genus $0$ as  $T_{0,n,\vec{r}} =
\G_{g,n,\vec{r}}$. The theorem is false in genus 1 as $\SL_2(\Z)$ does not have
the congruence subgroup property. In genus 2 the kernel of
$T_{2,n,\vec{r}}^\prol \to \G_{2,n,\vec{r}}^\prol$ is very large as we explain
in Section~\ref{sec:genus2}.

Theorem~\ref{thm:mcg_ladic} may be regarded as an integral analogue of results
about the relative unipotent completion of mapping class groups that are proved
in \cite{hain:comp,hain:torelli}. It is also an analogue for mapping class
groups of a result of Deligne \cite{deligne:central}, which may be stated as
follows. Suppose that $S$ is either $\SL_n$ ($n\ge 3$) or $\Sp_g$ ($g\ge 2$).
Suppose that $\G$ is a finite index subgroup of $S(\Z)$ and that $\Gtilde$ is
the restriction of the universal central extension of $S(\Z)$ to $\G$. Denote
the relative pro-$\ell$ completion of $\G$ (resp.\ $\Gtilde$) with respect to
the homomorphism to the closure of the image of $\G\to S(\Zl)$ (resp.\ $\Gtilde
\to S(\Zl)$) by $\G^\prol$ (resp.\ $\Gtilde^\prol$). Deligne's result implies
that there is an exact sequence
$$
0 \to \Zl \overset{\times e}{\longrightarrow} \Zl \to \Gtilde^\prol
\to \G^\prol \to 1
$$
where $e\in \{1,2\}$. It and the basic structure result for the relative
unipotent completion of mapping class groups \cite{hain:comp} mentioned above
are the two main tools in the proof of Theorem~\ref{thm:mcg_ladic}. The main
technical ingredient in the proof is the existence (Cor.~\ref{cor:cent_extn}) of
a central extension of a finite index subgroup of $\Sp_g(\Z)$ by $\Z$ in a
subquotient of $\G_{g,n,\vec{r}}$ whenever $g\ge 3$. This extension is a
non-zero multiple of the restriction of the universal central extension of
$\Sp_g(\Z)$.

Suppose that $n\ge 0$ and that $(S,\{x_0,\dots,x_n\})$ is an $n+1$ marked
surface of genus $g\ge 2$. Let $S'= S - \{x_1,\dots,x_n\}$ be the corresponding
$n$-punctured surface. Since the action of the mapping class group $\G_{g,n+1}$
on $\pi_1(S',x_0)$ induces an action of $\G_{g,n+1}^\prol$ on
$\pi_1(S',x_o)^\prol$, Theorem~\ref{thm:mcg_ladic} implies that the action of
the pro-$\ell$ Torelli group $T_{g,n+1}^\prol$ on $\pi_1(S',x_0)^\prol$ is not
faithful. Similarly, since the outer action of $\G_{g,n}$ on $\pi_1(S',x_0)$
induces an outer action of $\G_{g,n}^\prol$ on $\pi_1(S',x_0)^\prol$, the outer
action of $T_{g,n}^\prol$ on $\pi_1(S',x_0)^\prol$ is not
faithful.\footnote{When $g=2$, the statement follows from
Corollary~\ref{cor:ref}.}

\begin{bigcorollary}
If $g\ge 2$ and $n\ge 0$, then the natural homomorphisms
$$
T_{g,n+1}^\prol \to \Aut \pi_1(S',x_0)^\prol \text{ and }
T_{g,n}^\prol \to \Out \pi_1(S',x_0)^\prol
$$
are not injective. Their kernels each contain a copy of $\Zl$. \qed
\end{bigcorollary}

The injectivity of the corresponding homomorphisms
$$
\G_{g,n+1}^\prol \to \Aut \pi_1(S',x_0)^\prol \text{ and }
\G_{g,n}^\prol \to \Out \pi_1(S',x_0)^\prol
$$
for mapping class groups is not known for any $(g,n)$, with $g\ge 3$ and $n\ge
0$.

The proof of the following result (in the case $r=0$) was communicated to us by
Mochizuki and Tamagawa.

\begin{bigtheorem}
\label{thm:galois}
Suppose that $\ell$ is a prime number. For all non-negative integers $g,n,r$
satisfying $2g-2+n+r>0$, the action of $G_\Q$ on $\G_{g,n,\vec{r}}^\prol$,
determined by the choice of a suitably chosen smoothing of a maximally
degenerate stable curve of type $(g,n,r)$ as a tangential base point, is
unramified at all primes $\neq \ell$.
\end{bigtheorem}

We record the following consequence which follows directly from this result
using the functoriality of relative unipotent completion and the fact (proved in
\cite{hain:matsumoto_mcg} as a consequence of the ``comparison theorem'' proved
there) that the $G_\Q$-action on $\G_{g,n,\vec{r}}^\prol$ extends to the
relative unipotent completion $\cG_{g,n,\vec{r}}\otimes_\Q \Ql$ of
$\G_{g,n,\vec{r}}$ over $\Ql$.

\begin{bigcorollary}
For all prime numbers $\ell$, the $G_\Q$ representation on the relative
unipotent completion $\cG_{g,n,\vec{r}}\otimes \Ql$ is unramified outside
$\ell$. Consequently, the $G_\Q$ action on $\cG_{g,n,\vec{r}}\otimes \Ql$
induces an action of $\pi_1(\Spec\Z[1/\ell])$. \qed
\end{bigcorollary}

This result should be a useful technical ingredient in the investigation of
``Teichm\"uller motives'' --- relative motives over moduli spaces of curves ---
of which elliptic motives and elliptic polylogarithms \cite{beilinson-levin}
will be special cases.

\bigskip

\noindent{\em Acknowledgments:} The authors are grateful to S.~Mochizuki and
A.~Tamagawa for allowing us to reproduce the proof of Theorem~\ref{thm:galois}.
They are also grateful to Tamagawa for very helpful comments on an earlier draft
of this manuscript. Last, but not least, the authors owe an enormous debt of
gratitude the referee who read two drafts of the paper diligently. His/her
detailed and constructive comments and corrections have resulted in stronger
results, fewer errors and better exposition.

\section{Relative Pro-$\ell$ Completion}

\subsection{Basic Properties}

As in the introduction, we suppose that $\G$ is a discrete or profinite group,
that $P$ is a profinite group, and that $\rho : \G \to P$ is a continuous, dense
homomorphism.

The relative pro-$\ell$ completion of $\G$ with respect to $\rho$ is
characterized by a universal mapping property: if $G$ is a profinite group,
$\psi : G \to P$ a continuous homomorphism with pro-$\ell$ kernel, and if $\phi
: \G \to G$ is a continuous homomorphism whose composition with $\psi$ is
$\rho$, then there is a unique continuous homomorphism $\G^{\prol,\rho} \to G$
that extends $\phi$:
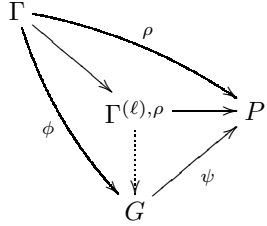
$$
\xymatrix{
\G \ar@/^/[drr]^\rho \ar@/_/[ddr]_\phi \ar[dr] \cr
& \G^{\prol,\rho} \ar[r]\ar@{.>}[d] & P \cr
& G \ar[ur]_\psi
}
$$

The following property is a direct consequence of the universal mapping
property.

\begin{proposition}
\label{prop:closure}
A dense homomorphism $\rho : \G \to P$ from a discrete group to a profinite
group induces a homomorphism $\rhobar : \Ghat \to P$ from the profinite
completion of $\G$ to $P$. The natural homomorphism $\G \to \Ghat$ induces a
natural isomorphism $\G^{\prol,\rho} \to \Ghat^{\prol,\rhobar}$. \qed
\end{proposition}

Suppose that $S$ is a linear group scheme over $\Z$. For a positive integer $m$
and a commutative ring $R$, set
$$
S(R)[m] = \ker\{S(R) \to S(R/mR)\}.
$$
This is the {\em level $m$ subgroup of $S(R)$}. Note that $S(\Zl)[m] =
S(\Zl)[\ell^\nu]$, where $\nu = \ord_\ell(m)$.

We say that $S(\Z)$ has the {\em congruence subgroup property} if every finite
index subgroup of $S(\Z)$ contains the level $m$ subgroup $S(\Z)[m]$ of $S(\Z)$
for some $m >0$, or equivalently, the natural homomorphism $S(\Z)^\wedge \to
S(\hat{\Z})$ is injective. In this case, the profinite completion of $S(\Z)$ is
the closure of its image in $S(\hat{\Z})$.

\begin{proposition}
\label{prop:S_ladic}
Suppose that $S$ is a linear group scheme over $\Z$. If $S$ satisfies
\begin{enumerate}

\item\label{item1}
$S(\Z)$ has the congruence subgroup property, and 

\item\label{item2}
$S(\Z)\to S(\Z/N)$ is surjective for every positive integer $N$,

\end{enumerate}
then the natural homomorphism $S(\Z)^\wedge \to S(\hat{\Z})$ is an isomorphism. 
If, in addition, $m$ is a positive integer and for every prime number $p$ not
dividing $\ell m$,  $S(\F_p)$ has no non-trivial quotient $\ell$-group, then,
the relative pro-$\l$ completion of $S(\Z)[m]$ with respect to $S(\Z)[m] \to
S(\Zl)[m]$ is the closure of $S(\Z)[m]$ in $S(\Zl)$, which equals
$S(\Zl)[\ell^\nu]$, where $\nu = \ord_\ell(m)$. 
\end{proposition}

\begin{proof}
The morphism
$$
S(\Z) \to S(\hat{\Z}) \cong \varprojlim_{N} S(\Z/N)
$$
induces a homomorphism $S(\Z)^\wedge \to S(\hat{\Z})$, which is injective
by condition (\ref{item1}) and surjective by condition (\ref{item2}).
Condition (\ref{item2}) also implies that the sequence
$$
1 \to {S(\Z)[m]} \to S(\Z) \to S(\Z/m\Z) \to 1
$$
is exact. Since the right group is finite, the sequence
$$
1 \to {S(\Z)[m]}^\wedge \to S(\Z)^\wedge \to S(\Z/m\Z) \to 1
$$
of profinite completions is exact. On the other hand, the sequence 
\begin{equation}
\label{eq:ses}
1 \to {S(\hat{\Z})[m]} \to S(\hat{\Z}) \to S(\Z/m\Z) \to 1
\end{equation}
is exact by definition. The isomorphism $S(\Z)^\wedge \to S(\hat{\Z})$ therefore
induces an isomorphism  of $S(\Z)[m]^\wedge$ with $S(\hat{\Z})[m]$.

Denote the order of $m$ at a prime number $p$ by $\nu_p$.  The short exact
sequence (\ref{eq:ses}) is the direct product over all prime numbers $p$ of the
exact sequences
$$
1 \to {S(\Z_p)[p^{\nu_p}]} \to S(\Z_p) \to S(\Z/p^{\nu_p}\Z) \to 1.
$$
Consequently, the profinite completion of $S(\Z)[m]$ is
$$
S(\hat{\Z})[m] = \prod_p S(\Z_p)[p^{\nu_p}]
$$
where $p$ ranges over all prime numbers. The homomorphism $\rho : S(\Z)[m] \to
S(\Zl)[m]$ induces the projection
$$
S(\hat{\Z})[m] \cong \prod_p S(\Z_p)[p^{\nu_p}] \to S(\Zl)[\ell^{\nu_\ell}]
$$
onto the factor corresponding to $\l$. Since the relative pro-$\ell$ completion 
$S(\Z)[m]^\prol$ is the largest quotient $G$ of the profinite completion
$S(\hat{\Z})[m]$ such that the kernel of  $G \to S(\Zl)[m]$ is pro-$\ell$, there
is a surjection
$$
\prod_{p \neq \ell} S(\Z_p)[p^{\nu_p}]
\to \ker\{S(\Z)[m]^\prol \to S(\Zl)[m]\} 
$$
whose image is pro-$\ell$. Thus, to prove the proposition, it suffices to show
that $S(\Z_p)[p^{\nu_p}]$ has no non-trivial $\ell$-quotient when $p \neq
\ell$. 

Now suppose that $p\neq \ell$. It is well known that $S(\Z_p)[p^N]$ is a pro-$p$
group when $N>0$. (It suffices to check this for $\GL_M$.) So if $N>0$, every
homomorphism from $S(\Z_p)[p^N]$ to an $\ell$-group is trivial. In particular,
if $\nu_p > 0$ (or equivalently if $p$ divides $m$), then $S(\Z_p)[p^{\nu_p}]$ 
has no non-trivial $\ell$-quotients, as required.  Suppose that $p$ does not
divide $m\l$. In particular, $\nu_p=0$. Since
$$
1 \to S(\Z_p)[p] \to S(\Z_p) \to S(\F_p) \to 1
$$
is exact by condition (\ref{item2}), every homomorphism $S(\Z_p) \to A$ to an
$\ell$-group ($\ell\neq p$) factors through a homomorphism $S(\F_p) \to A$. But
we have assumed that there are no non-trivial such homomorphisms. Thus
$S(\Z_p)[p^{\nu_p}]$ has no non-trivial $\ell$-quotients, as required.
\end{proof}

The hypotheses of Proposition~\ref{prop:S_ladic} are satisfied for all $m\ge 1$
when $S=\Sp_g$ when $g \geq 2$ and when $S=\SL_n$ when $n\ge 3$. (See \cite{bms}
and \cite[Thm.~1, p.~12]{dieudonne}.)

Define $\Gp_f$ to be the category whose objects are topological groups  that are
either discrete or profinite; morphisms are continuous group homomorphisms. The
category of discrete groups and the category of profinite groups are both full
subcategories of $\Gp_f$. A sequence of group homomorphisms in $\Gp_f$ is exact
if its is exact in the category of groups.

\begin{proposition}[naturality]
Suppose that $\G_1$ and $\G_2$ are objects of $\Gp_f$ and that $P_1$ and $P_2$
are profinite groups. Suppose that $\rho_j : \G_j \to P_j$ are continuous dense
homomorphisms. If
$$
\xymatrix{
\G_1 \ar[d]_{\phi_\G} \ar[r]^{\rho_1} & P_1 \ar[d]^{\phi_P}\cr
\G_2 \ar[r]^{\rho_2} & P_2 
}
$$
is a commutative diagram in $\Gp_f$, then there is a unique continuous
homomorphism $\phi^\prol : \G_1^{\prol,\rho_1} \to \G_2^{\prol,\rho_2}$ such
that the diagram
$$
\xymatrix{
\G_1 \ar[d]_{\phi_\G} \ar[r] \ar@/^1pc/[rr]^{\rho_1} &
\G_1^{\prol,\rho_1}\ar[d]_{\phi^\prol} \ar[r] & P_1 \ar[d]_{\phi_P}\cr
\G_2 \ar@/_1pc/[rr]_{\rho_2} \ar[r] & \G_2^{\prol,\rho_2}\ar[r] & P_2 
}
$$
commutes.
\end{proposition}

\begin{proof}
This follows easily from the universal mapping property.
\end{proof}

\begin{proposition}[right exactness]
\label{prop:right-exact-prol}
Suppose that $P_1$, $P_2$ and $P_3$ are profinite groups and that $\rho_j : \G_j
\to P_j$ ($j \in \{1,2,3\}$) are continuous dense homomorphisms in $\Gp_f$.
Suppose that the $\G_j$ are all discrete groups or all profinite groups. If the
diagram
$$
\xymatrix{
1 \ar[r] & \G_1 \ar[r]\ar[d]_{\rho_1} &
\G_2 \ar[r]\ar[d]_{\rho_2} & \G_3 \ar[r]\ar[d]_{\rho_3} & 1 \cr
1 \ar[r] & P_1 \ar[r] & P_2 \ar[r] & P_3 \ar[r] & 1
}
$$
in $\Gp_f$ commutes and has exact rows, then the sequence
$$
\G_1^{\prol,\rho_1} \to \G_2^{\prol,\rho_2} \to \G_3^{\prol,\rho_3} \to 1
$$
is exact.
\end{proposition}

\begin{proof}
Denote $\G_j^{\prol,\rho_j}$ by $\G_j^\prol$. Since $\G_2 \to \G_3$ is
surjective, the homomorphism $\G_2 \to \G_3 \to \G_3^\prol$ has dense image.
This implies that $\G_2^\prol \to \G_3^\prol$ is surjective.

Denote the closure of the image of $\G_1 \to \G_2^\prol$ by $N$. Observe that
this equals the image of $\G_1^\prol \to \G_2^\prol$. Since $\G_1$ is normal in
$\G_2$, this is a closed normal subgroup of $\G_2^\prol$. The image of $N$ in
$P_2$ is $P_1$. The homomorphism $\G_2^\prol \to P_2$ induces a homomorphism
$\rhobar : \G_2^\prol/N \to P_2/P_1 \cong P_3$.  Since the kernel of the
composite $\G_2 \to \G_2^\prol \to \G_2^\prol/N$ contains $\G_1$ it induces a
continuous homomorphism $\rho : \G_3 \to \G_2^\prol/N$ whose composition with
$\rhobar$ is $\rho_3$. The universal mapping property of $\G_3 \to \G_3^\prol$
implies that $\rho$ induces a homomorphism $\G_3^\prol \to \G_2^\prol/N$, which
is necessarily an isomorphism. Since $N$ is the image of $\G_1^\prol \to
\G_2^\prol$, this implies the exactness at $\G_2^\prol$.
\end{proof}

Relative pro-$\l$ completion is not left exact. One dramatic illustration of
this is Deligne's result \cite{deligne:central}, which we now describe.

Suppose that $S$ is the $\Z$-group scheme $\SL_n$ or $\Sp_g$ where $n\ge 3$ or
$g\ge 2$. Since $S(\Z)$ is perfect, $S(\Z)$ has a universal central extension.
Denote the restriction of the universal central extension of $S(\Z)$ to a finite
index subgroup $\G$ by $\Gtilde$. This is an extension
$$
0 \to \Z \to \Gtilde \to \G \to 1.
$$
Deligne \cite{deligne:central} proves that every finite index subgroup of
$\Gtilde$ contains $2\Z$, which implies that the profinite completion of
$\Gtilde$ is an extension
$$
0 \to \Z/e\Z \to \Gtilde^\wedge \to \G^\wedge \to 1
$$
where $e\in \{1,2\}$. Since $S(\Z)$ has the congruence subgroup property
\cite{bms}, the profinite completion of $\G$ is its closure in $S(\hat{\Z})$.

Deligne's result can be reformulated and extended to multiples of the universal
central extension $\Gtilde \to \G$.

Suppose as above that $S = \SL_n$ or $\Sp_g$ with $n\ge 3$ or $g\ge 2$.  Suppose
that $\G$ is a finite index subgroup of $S(\Z)$. Denote the restriction to $\G$
of the $d$th power of the universal central extension of $S(\Z)$ by $\Gtilde_d$.
(When $d=0$, $\Gtilde_d = \G\times \Z$.) Note that $\Gtilde = \Gtilde_1$. For
all non-zero $d\in \Z$ there is a short exact sequence
$$
1 \to \Gtilde \to \Gtilde_d \to \Z/d\Z \to 0.
$$

When the context is clear, we will often drop $\rho$ from the notation and
denote $\G^{\prol,\rho}$ by $\G^\prol$. For example, we will denote the 
relative pro-$\l$ completion of $\Sp_n(\Z)$ with respect to reduction mod $\l$
by $\Sp_n(\Z)^\prol$.

\begin{proposition}
\label{prop:cent-extn}
For all prime numbers $\ell$ and all non-zero $d\in \Z$, there are short exact
sequences
$$
1 \to \Gtilde^\prol \to \Gtilde_d^\prol \to \Zl/d\Zl \to 0
$$
and
$$
0 \to \Zl/df\Zl \to \Gtilde_d^\prol \to \G^\prol \to 1,
$$
where $f\in\{1,2\}$ and $f\le e$, the constant defined above.
\end{proposition}

\begin{proof}
Applying right exactness to the standard exact sequence
$$
1 \to \Gtilde \to \Gtilde_d \to \Z/d\Z \to 0
$$
implies that the sequence
$$
\Gtilde^\prol \to \Gtilde_d^\prol \to \Zl/d\Zl \to 0
$$
is exact. Denote the kernel of $\Gtilde_d \to \G$ by $K_d$.  Set $K = K_1$. Then
each $K_d$ is isomorphic to $\Z$ and $K_1$ is the subgroup of $K_d$ of index
$d$. Injectivity of $\Gtilde^\prol \to \Gtilde_d^\prol$ follows from the fact
that if $N$ is a normal subgroup of $\Gtilde$ whose intersection with $K$ has
$\l$-power index,  then the subgroup $\tilde{N}:=N\cdot f^{-1}(K \cap N)$ of
$\Gtilde_d$, where $d = \ell^\nu f$ with $\nu = \ord_\l(d)$, is normal in
$\Gtilde_d$ and has the property that its intersection with $K_d$ has $\l$-power
index and $\tilde{N}\cap \Gtilde = N$.

The diagram
$$
\xymatrix{
& \Z/d\Z \ar@{=}[r] & \Z/d\Z \cr
0 \ar[r] & \Z \ar[r]\ar@{->>}[u] & \Gtilde_d \ar[r]\ar@{->>}[u]
& \G\ar@{=}[d] \ar[r]& 1 \cr
0 \ar[r] & \Z \ar[r]\ar[u]_{\times d} & \Gtilde \ar[r]\ar[u] & \G \ar[r] & 1
}
$$
has exact rows; the upper vertical maps are surjective and the lower vertical
maps are injective. Deligne's result and the first part of this result imply
that the relative pro-$\l$ completion of this diagram is
$$
\xymatrix{
& \Zl/d\Zl \ar@{=}[r] & \Zl/d\Zl \cr
 & \Zl/df\Zl \ar[r]\ar@{->>}[u] & \Gtilde_d^\prol \ar[r]\ar@{->>}[u]
& \G^\prol\ar@{=}[d] \ar[r]& 1 \cr
0 \ar[r] & \Zl/f\Zl \ar[r]\ar[u]_{\times d} & \Gtilde^\prol \ar[r]\ar[u]
& \G^\prol \ar[r] & 1
}
$$
where $f\le e$ and $f\in \{1,2\}$. The second assertion follows.
\end{proof}

The following useful lemma is easily proved using the universal mapping property
of relative pro-$\ell$ completion.

\begin{lemma}
Suppose that
$$
\xymatrix{
1 \ar[r] & K \ar[r] & P \ar[r]^{\psi} & \overline{P} \ar[r] & 1
}
$$
is a short exact sequence of profinite groups. Suppose that $\G$ is a discrete
or profinite group and that $\rho : \G \to P$ is a continuous dense
homomorphism. Denote $\psi\circ \rho$ by $\rhobar$. If $K$ is a pro-$\ell$
group, then the natural homomorphism $\G^{\prol,\rho} \to \G^{\prol,\rhobar}$ is
an isomorphism. \qed
\end{lemma}

\begin{example}
\label{ex:complns_equal}
Let $S$ be the $\Z$-group scheme $\SL_{n+1}$ or $\Sp_n$, where $n\ge 1$. Let
$P=S(\Zl)$ and $\overline{P} = S(\Fl)$. Suppose that $\G$ is a discrete or
profinite group and that $\rho : \G \to S(\Zl)$ is a continuous, dense
homomorphism. Suppose that $\psi : P \to \overline{P}$ is reduction mod $\l$.
Since $\ker \psi$ is a  pro-$\l$ group, there is a natural isomorphism
$\G^{\prol,\rho} \cong \G^{\prol,\rhobar}$.
\end{example}

\subsection{Remarks on pro-$\l$ completion}
\label{sec:ladic}
To conclude this section, we establish some basic facts about standard (i.e.,
relative to the trivial representation) pro-$\l$ completion. Denote the lower
central series of a discrete group $\G$ by
$$
\G = L^1 \G \supseteq L^2\G \supseteq L^3\G \supseteq \cdots
$$

The following observation follows directly from the fact that every finite
$\ell$-group is nilpotent.

\begin{lemma}
\label{lem:prol}
For all discrete groups $\G$, the natural mapping
$$
\G^\prol \to \varprojlim_n(\G/L^n)^\prol
$$
is an isomorphism. \qed
\end{lemma}

This result can be refined as follows. For an integer $m$, let $\supp(m)$ denote
the set of prime divisors of $m$. For a set $S$ of prime numbers and a subgroup
$N$ of $\G$, define
$$
\sqrt[S]{N} = \{g \in \G : \text{ there exists an integer $m>0$ with }
\supp(m) \subseteq S \text{ and } g^m  \in N\}.
$$
The fact that the set of torsion elements of a nilpotent group is a
characteristic subgroup, implies that if $N$ is a normal (resp.\ characteristic)
subgroup of $\G$ and $\G/N$ is nilpotent, then for all $S$, $\sqrt[S]{N}$ is a
normal (resp.\ characteristic) subgroup of $\G$.

Denote the complement of $S$ in the set $\wp$ of prime numbers by $S'$. When $S$
consists of a single prime number $\l$, we shall write $\sqrt[\l]{N}$ for
$\sqrt[S]{N}$ and $\sqrt[\l']{N}$ for $\sqrt[S']{N}$. For simplicity, we denote
$\sqrt[\wp]{N}$ by $\sqrt{N}$. This is the set of all elements of $\G$ that are
torsion mod $N$.

Again, the fact that every finite $\l$ group is nilpotent implies:

\begin{lemma}
\label{lem:prol_bis}
For all discrete groups $\G$, the natural mapping
$$
\G^\prol \to \varprojlim_n(\G/\sqrt[\l']{L^n})^\prol
$$
is an isomorphism. \qed
\end{lemma}

For each $n \ge 1$, define $D^n \G$ to be the inverse image in $\G$ of the
torsion elements of $\G/L^n \G$. That is, $D^n\G = \sqrt{L^n \G}$, which is a
characteristic subgroup of $\G$. It is known as the $n$th {\em rational
dimensional subgroup of $\G$}.\footnote{The groups $D^n\G$ can also be
constructed as follows. Suppose that $F$ is a field of characteristic zero.
Denote the augmentation ideal of the group algebra $F\G$ by $J_F$. Then $D^n\G$
is the intersection of $\G$ with $1+J_F^n$ in $F\G$. For a proof see
\cite[Thm.~1.10, p.~474]{passman}.}

Recall from the introduction that a discrete group $\G$ is $\l$-regular if every
finite quotient of $\G$ that is an $\l$-group is the quotient of a torsion free
nilpotent quotient of $\G$.

The condition that a discrete group $\G$ be $\l$-regular may be expressed in
terms of the groups $D^m\G$ and $\sqrt[\l']{L^n\G}$. Observe that for all
positive integers $n$, $D^n\G \supseteq \sqrt[\l']{L^n\G}$.

\begin{proposition}
\label{prop:comp_rtfn}
A discrete group $\G$ is $\l$-regular if for each positive integer $n$, there
exists an integer $m>0$ such that $D^m\G \subseteq \sqrt[\l']{L^n\G}$. In this
case, the mapping
$$
\G^\prol \to \varprojlim_n (\G/D^n)^\prol
$$
is an isomorphism.
\end{proposition}

\begin{proof}
The natural mapping
$$
\G^\prol \cong \varprojlim_n(\G/\sqrt[\l']{L^n})^\prol
\to \varprojlim_n (\G/D^n)^\prol
$$
is surjective. The condition in the statement guarantees that it is also
injective.
\end{proof}

Not all groups that inject into their pro-$\l$ completions are $\l$-regular.

\begin{example}
Suppose that $S = \SL_{n+1}$ or $\Sp_n$ where $n\ge 2$. Suppose that $\G =
S(\Z)[\l^m]$ where $m\ge 1$. Even though this group injects into its pro-$\l$
completion $S(\Zl)[\l^m]$, the homomorphism
$$
\G^\prol \to \varprojlim_k (\G/D^k)^\prol
$$
is trivial as $\G$ has no non-trivial torsion free abelian quotients
\cite{raghunathan}, and therefore no non-trivial torsion free nilpotent
quotients.
\end{example}

If the graded quotients of the lower central series of $\G$ have no
$\l$-torsion, then $D^n\G = \sqrt[\l']{L^n\G}$ for all $n>0$. Thus we have:

\begin{corollary}
\label{cor:comp_rtfn}
Suppose that $\G$ is a discrete group. If the graded quotients of the lower
central series of $\G$ have no $\l$-torsion, then $\G$ is $\l$-regular. \qed
\end{corollary}

We conclude this section with the following result, which follows easily from
the fact that every finite $\l$-group is nilpotent.

\begin{proposition}
If $\G$ is a discrete group such that $H_1(\G)\otimes \Fl = 0$, then $\G^\prol$
is trivial. \qed
\end{proposition}

\section{Mapping Class and Torelli Groups}

Suppose that $g$ and $n$ are non-negative integers satisfying $2g-2+n > 0$. Fix
a closed, oriented surface $S$ of genus $g$ and a finite subset $\x =
\{x_1,\dots,x_n\}$ of $n$ distinct points in $S$. The corresponding mapping
class group will be denoted
$$
\G_{S,\x} = \pi_0\Diff^+ (S,\x),
$$
where $\Diff^+(S,\x)$ denotes the group of orientation preserving
diffeomorphisms of $S$ that fix $\x$ pointwise. By the classification of
surfaces, the diffeomorphism class of $(S,\x)$ depends only on $(g,n)$.
Consequently, the group $\G_{S,\x}$ depends only on the pair $(g,n)$. It will be
denoted by $\G_{g,n}$ when we have no particular marked surface $(S,\x)$ in
mind.

For a commutative ring $A$, set $H_A = H_1(S;A)$. The intersection pairing
$H_A^{\otimes 2} \to A$ is skew symmetric and unimodular. The choice of a
symplectic basis of $H_A$ gives an isomorphism $\Sp(H_A) \cong \Sp_g(A)$ of
$2g\times 2g$ symplectic matrices with entries in $A$. The action of
$\Diff^+(S,\x)$ on $S$ induces a homomorphism
$$
\rho : \G_{S,\x} \to \Sp(H_\Z)
$$ 
which is well-known to be surjective.

Let $\rhobar : \G_{S,\x} \to \Sp(H_{\Z/\l\Z})$ be reduction of $\rho$ mod $\l$.
Define $\G_{S,\x}^\prol$ to be the relative pro-$\l$ completion of $\G_{S,\x}$
with respect to $\rhobar$. By Example~\ref{ex:complns_equal}, this is isomorphic
to the relative pro-$\l$ completion of $\G_{S,\x}$ with respect to the natural
homomorphism $\G_{S,\x} \to \Sp(H_\Zl)$.

\begin{proposition}
\label{prop:injectivity}
Suppose that $n\ge 0$, that $S$ is a surface of genus $\ge 2$, and that $\x =
\{x_0,\dots,x_n\}$ is a set of $n+1$ distinct points of $S$. Set $\x' = \x
-\{x_0\}$ and $S' = S - \x'$. Then
\begin{enumerate}

\item the natural homomorphism $\G_{S,\x} \to \G_{S,\x}^\prol$ is injective;

\item the sequence $1 \to \pi_1(S',x_0)^\prol \to \G_{S,\x}^\prol \to
\G_{S,\x'}^\prol \to 1$ of relative pro-$\l$ completions is exact.

\end{enumerate}
\end{proposition}

\begin{proof}
Since $\pi_1(S',x_0)$ is residually nilpotent and since the graded quotients of
its lower central series are torsion free (cf.\ \cite{labute} when $n=0$), the
homomorphism $\pi_1(S',x_0) \to \pi_1(S',x_0)^\prol$ is injective. Thus
$\Aut \pi_1(S',x_0) \to \Aut_\cts \pi_1(S',x_0)^\prol$ is also injective. Since
$\G_{S,\x}$ is a subgroup of $\Aut \pi_1(S',x_0)$, it is a subgroup of
$\Aut_\cts \pi_1(S',x_0)^\prol$.  The first assertion now follows as the
inclusion of $\G_{S,\x}$ into $\Aut_\cts \pi_1(S',x_0)^\prol$ factors through
the completion homomorphism:
$$
\G_{S,\x} \to \G_{S,\x}^\prol \to \Aut_\cts \pi_1(S',x_0)^\prol.
$$
Since relative pro-$\l$ completion is right exact, to prove the second
assertion, we need only show that $\pi_1(S',x_0)^\prol \to \G_{S,\x}^\prol$ is
injective. But this follows as the composite
$$
\pi_1(S',x_0)^\prol \to \G_{S,\x}^\prol \to \Aut_\cts \pi_1(S',x_0)^\prol
$$
is the conjugation action, which is injective as $\pi_1(S',x_0)^\prol$ has
trivial center (cf.\ \cite{anderson} when $n=0$).
\end{proof}

The {\em Torelli group} $T_{S,\x}$ is defined to be the kernel of $\rho$. Its
isomorphism class depends only on $(g,n)$. It will be denoted by $T_{g,n}$ when
we have no particular marked surface in mind.

\begin{proposition}
Suppose that $n\ge 0$, that $S$ is a surface of genus $\ge 2$, and that $\x =
\{x_0,\dots,x_n\}$ is a set of $n+1$ distinct points of $S$. Set $\x' = \x
-\{x_0\}$ and $S' = S - \x'$. Then the sequence
$$
1 \to \pi_1(S',x_0)^\prol \to T_{S,\x}^\prol \to T_{S,\x'}^\prol \to 1
$$
of relative pro-$\l$ completions is exact.
\end{proposition}

\begin{proof}
We need only prove the injectivity of the left hand mapping. But this follows
from the injectivity of the composite
$$
\pi_1(S',x_0)^\prol \to T_{S,\x}^\prol \to \G_{S,\x}^\prol
$$
which was established in Proposition~\ref{prop:injectivity}.
\end{proof}

\subsection{Variant: tangent vectors and boundary components}
\label{sec:mcgs}
Suppose that $\x$ is a set of $n$ distinct points in $S$ and that $v_j$ is a
non-zero tangent vector at $y_j \in S$, and $\{y_1,\dots,y_r\}$ is a set of $r$
distinct points, disjoint from $\x$. Set $\vecs = \{v_1,\dots,v_r\}$. Define
$$
\G_{S,\x,\vecs} = \pi_0 \Diff^+(S,\x,\vecs),
$$
the group of connected components of the group of orientation preserving
diffeomorphisms of $S$ that fix $\x$ and $\vecs$ pointwise. The corresponding
Torelli group $T_{S,\x,\vecs}$ is the kernel of the homomorphism
$\G_{S,\x,\vecs} \to \Sp(H_1(S))$.  For $m\in \N$, define its level $m$ subgroup
$\G_{g,n,\vec{r}}[m]$ to be the kernel of the natural homomorphism
$\G_{g,n,\vec{r}} \to \Sp(H_1(S;\Z/m))$.

The group $\G_{S,\x,\vecs}$ is isomorphic to the mapping class group of a genus
$g$ surface with $n$ marked points and $r$ boundary components. The isomorphism
can be seen by replacing $S-\{x_1,\dots,x_n,y_1,\dots,y_r\}$ by the real
oriented blow-up of $S-\{x_1,\dots,x_n\}$ at each of the $y_j$. (See \cite[\S
4.1]{hain:morita} for more details.)

The groups $\G_{S,\x,\vecs}$ and $T_{S,\x,\vecs}$ depend only on $g$, $n$ and
$r$. We shall often denote them by $\G_{g,n,\vec{r}}$ and $T_{g,n,\vec{r}}$. The
indices $n$ and $\vec{r}$ will be dropped when they are zero. For example,
$\G_{g,\vec{1}}$ denotes the mapping class group associated to a genus $g$
surface with one tangent vector, while $\G_{g,1}$ denotes the mapping class
group of a genus $g$ surface with one marked point.

Replacing each tangent vector by its anchor point defines a natural surjective
homomorphism $\G_{g,n,\vec{r}} \to \G_{g,n+r}$ whose kernel consists of the
free abelian group generated by the $r$ Dehn twists about the anchor points
of the tangent vectors. The extension
$$
0 \to \Z^r \to \G_{g,n,\vec{r}} \to \G_{g,n+r} \to 1
$$
is central and, by right exactness of relative pro-$\l$ completion, induces an
exact sequence
$$
\Zl^r \to \G_{g,n,\vec{r}}^\prol \to \G_{g,n+r}^\prol \to 1.
$$

\begin{proposition}
\label{prop:injective}
The sequence
$$
0 \to \Zl^r \to \G_{g,n,\vec{r}}^\prol \to \G_{g,n+r}^\prol \to 1
$$
is exact.
\end{proposition}

\begin{proof}
Suppose that $2g-1+m>0$. Denote the fundamental group of a genus $g$ surface $S$
with $m$ punctures and one boundary component (with base point lying on the
boundary) by $\Pi_{g,m,\vec{1}}$. It is a free group $\langle
a_1,\dots,a_g,b_1,\dots,b_g,z_1,\dots,z_m\rangle$ of rank $2g+m$ generated by a
standard set of generators. The boundary loop $c$ is given by
$$
c = z_1^{-1}\cdots z_m^{-1}\prod_{j=1}^g a_j b_j a_j^{-1}b_j^{-1}.
$$
The action of the Dehn twist on the boundary is conjugation by $c^{-1}$, which
is of infinite order as $\Pi_{g,m,\vec{1}}$ is free.

Since $\Pi_{g,m,\vec{1}}$ is free, it injects into its pro-$\l$ completion
$\Pi_{g,m,\vec{1}}^\prol$. Consequently, the action of a Dehn twist about the
boundary of $S$ on $\Pi_{g,m,\vec{1}}^\prol$ has infinite order. In other words,
the restriction of the natural homomorphism
$$
\psi : \G_{g,m,\vec{1}} \to \Aut \Pi_{g,m,\vec{1}}^\prol
$$
to the central $\Z$ generated by the Dehn twist about the boundary, is
injective. This proves the result when $r=1$.

When $r>1$, consider the homomorphism $\phi_j : \G_{g,n,\vec{r}} \to
\G_{g,n+r-1,\vec{1}}$ that replaces all but the $j$th ($1\le j \le r$) tangent
vector by its anchor point. The composite
$$
\xymatrix{
\Z^r \ar[r] & \G_{g,n,\vec{r}} \ar[r]^(0.35){\psi\circ \phi_j} &
\Aut \Pi_{g,n+r-1,\vec{1}}^\prol
}
$$
is injective on the $j$th factor, and trivial on all others. This implies
that the homomorphism $\Zl^r \to \G_{g,n,\vec{r}}^\prol$ is injective.
\end{proof}

\subsection{Continuous actions}

For later use, we show that when $g\ge 3$, even though $T_{g,n,\vec{r}}^\prol$
may not inject into $\G_{g,n,\vec{r}}^\prol$, there is still a natural action of
$\G_{g,n,\vec{r}}^\prol$ on $T_{g,n,\vec{r}}^\prol$.

\begin{proposition}
\label{prop:induced_action}
If $g\ge 3$, then for all $r,n \ge 0$, the conjugation action $\G_{g,n,\vec{r}}
\to \Aut T_{g,n,\vec{r}}$ induces a continuous action $\G_{g,n,\vec{r}}^\prol
\to \Aut T_{g,n,\vec{r}}^\prol$.
\end{proposition}

\begin{proof}
Set $\G = \G_{g,n,\vec{r}}$, $T = T_{g,n,\vec{r}}$ and $S=\Sp_g$. Denote the
genus $g$ reference surface by $X$. Since pro-$\ell$ completion is functorial,
the conjugation action induces an action $\G \to \Aut T^\prol$. Since $T$ acts
trivially on the graded quotients of the lower central series of $T^\prol$, its
image in $\Aut T^\prol$ lies in the kernel $\Aut_0 T^\prol$ of the natural
homomorphism $\Aut T^\prol \to \Aut H_1(T^\prol)$. Since $g\ge 3$, $T$ is
finitely generated by \cite{johnson:fg}, which implies that $H_1(T^\prol) =
H_1(T)\otimes\Zl$. Since the sequence
$$
1 \to \Aut_0 T^\prol \to \Aut T^\prol \to \Aut H_1(T)\otimes \Zl
$$
is exact and $\Aut_0 T^\prol$ is a pro-$\ell$ group, to prove the result, it
suffices to show that the action of $S(\Z)$ on $H_1(T)\otimes\Zl$ factors
through $S(\Z) \to S(\Zl)$. Johnson \cite{johnson:h1} has shown that $H_1(T)$
has only 2-torsion and that $H_1(T)/\text{torsion}$ is an algebraic
representation of the group scheme $S$. This implies that when $\ell$ is odd,
the $S(\Z)$ action on $H_1(T)\otimes \Zl$ factors through $S(\Z) \to S(\Zl)$.

To establish the result when $\ell = 2$, we need another fundamental fact that
was established by Johnson in \cite{johnson:h1}. Namely, $H_1(T;\F_2)$ is
isomorphic to a space\footnote{The space depends on $n$ and $r$.} of cubic
boolean polynomials on the $\F_2$-valued quadratic forms on $H_1(X;\F_2)$ with
respect to its intersection form. This implies that the action of $S(\Z)$ on
$H_1(T;\F_2)$ factors through $S(\Z) \to S(\F_2)$. Since the $2$-torsion
$H_1(T)_{(2)}$ of $H_1(T)$ is an $S(\Z)$ submodule of $H_1(T;\F_2)$, it is also
an $S(\F_2)$-module. To show that the action of $S(\Z)$ on $H_1(T)\otimes \Z_2$
factors through $S(\Z) \to S(\Z_2)$, it suffices to show that the image of
$S(\Z)[2]$ in $\Aut H_1(T)\otimes \Z_2$ is contained in a pro-$2$ subgroup of
$\Aut H_1(T)\otimes \Z_2$. This is because $S(\Z_2)$ is the relative pro-$2$
completion  of $S(\Z)$ with respect to $S(\Z) \to S(\F_2)$ by 
Proposition~\ref{prop:S_ladic} and Example~\ref{ex:complns_equal}.

Write $H_1(T)$ as an extension
$$
0 \to H_1(T)_{(2)} \to H_1(T) \to V \to 0.
$$
Johnson's computation of $H_1(T)$ implies that $V$ is a finitely generated,
torsion free $\Z$-module that is an algebraic representation of $S$. The image
of $S(\Z)[2]$ in $\Aut H_1(T)$ preserves this sequence and is thus contained in
a subgroup $G$ of $\Aut H_1(T)\otimes \Z_2$ that is an extension
$$
0 \to \Hom_\Z(V,H_1(T)_{(2)}) \to G \to S(\Z_2)[2] \to 1,
$$
which is a pro-$2$ group and thus completes the proof.
\end{proof}

\subsection{Notational Convention}

To avoid confusion, we make explicit our convention that the {\em relative}
pro-$\l$ completion of the mapping class group $\G_{g,n,\vec{r}}[m]$ (resp.\ the
arithmetic group $\Sp_g(\Z)[m]$) with respect to the natural homomorphism to
the closure of its image in $\Sp_g(\Zl)$ will be denoted by
$\G_{g,n,\vec{r}}[m]^\prol$ (resp.\ $\Sp_g(\Z)[m]^\prol$).

\section{Relative Unipotent Completion of Discrete Groups}

Relative unipotent completion is analogous to relative pro-$\l$ completion, but
often more computable as it is controlled by cohomology. In many situations,
including the ones in this paper, the relative pro-$\l$ completion of a group
$\G$ maps to its relative unipotent completion. In such cases, relative
unipotent completion can be used to give a lower bound for the size of the
relative pro-$\l$ completion of a discrete group; under very favorable
conditions, it can be used to compute the relative pro-$\l$ completion.

There are various kinds of relative unipotent completion, in this paper we need
only the simplest kind, which is reviewed in this section. More details can be
found in \cite{hain:comp,hain:morita}. Throughout this section, $F$ will denote
a field of characteristic zero.

The basic data are:
\begin{enumerate}

\item a discrete group $\G$;

\item a reductive linear algebraic $F$-group $R$;

\item a Zariski dense representation $\rho : \G \to R(F)$.

\end{enumerate}
The relative unipotent completion of $\G$ with respect to $\rho$ is a
proalgebraic group $\cG$ over $F$ which is an extension
$$
1 \to \U \to \cG \to R \to 1,
$$
where $\U$ is a prounipotent $F$-group, together with a homomorphism $\rhohat :
\G \to \cG(F)$ that lifts $\rho$. It is characterized by the following universal
mapping property: if $G$ is a linear algebraic $F$-group that is an extension of
$R$ by a unipotent group and if $\rhotilde : \G \to G(F)$ is a homomorphism that
lifts $\rho$, then there is a unique homomorphism $\phi : \cG \to G$ of
proalgebraic $F$-groups such that
$$
\xymatrix@!C{
& \cG(F) \ar[dr] \ar[dd]_\phi \cr
\G \ar[ur]^\rhohat \ar[dr]_\rhotilde && R(F) \cr
& G(F) \ar[ur]
}
$$
commutes.

We shall denote the relative unipotent completion of $\G$ with respect to $\rho$
by $\cG(\G,\rho)$, or simply $\cG(\G)$ when $\rho$ is clear from context.

The relative unipotent completion of a discrete group can be constructed as an
inverse limit of Zariski dense representations $\rhohat : \G \to G(F)$, where
the $F$-group $G$ is an extension of $R$ by a unipotent group. Alternatively,
$\cG(\G)$ can be constructed as the tannakian fundamental group of the category
of finite dimensional $F\G$-modules $V$ that admit a filtration
$$
0 = I_0 V \subseteq I_1 V \subseteq  \cdots \subseteq I_{n-1}V \subseteq
I_n V = V 
$$
by $F\G$-submodules with the property that each graded quotient $I_kV/I_{k-1}V$
has a rational action of $R$ for which the action of $\G$ on it factorizes
through $\rho : \G \to R(F)$. This is a neutral tannakian category
$\cR(\G,\rho)$ over $F$. There is a natural isomorphism
$$
\pi_1(\cR(\G,\rho),\omega) \cong \cG(\G,\rho)
$$
where $\omega : \cR(\G,\rho) \to \Vect_F$ is the forgetful functor. See
\cite{deligne-milne} for definitions and \cite[\S 7]{hain-matsumoto:survey} for
discussion of tannakian description of weighted completion.

\begin{example}
\label{ex:lattice}
If $\G$ is an arithmetic subgroup of an almost simple $\Q$-group $G$ of real
rank $\ge 2$, and $F$ is any field of characteristic zero, then the relative
unipotent completion of $\G$ with respect to the inclusion $\G \hookrightarrow
G(F)$ is $G_{/F}$, \cite[p.~84]{hain:comp}. In particular, if $\G$ is a finite
index subgroup of $\Sp_g(\Z)$ and $g\ge 2$ (resp.\ of $\SL_n(\Z)$ and $n\ge 3$),
then the relative unipotent completion of $\G$ with respect to the inclusion
$\rho: \G \to \Sp_g(F)$ (resp.\ $\G \to \SL_n(F)$) is $\Sp_{g/F}$ (resp.\
$SL_{n/F}$). When $g=1$, the completion of $\G$ has a very large prounipotent
radical --- it is free of countable rank, \cite[Remark~3.9]{hain:torelli}.

\subsection{Unipotent completion}
When $R$ is the trivial group, relative unipotent completion reduces to {\em
unipotent (also called Malcev) completion}. We shall denote the unipotent
completion of $\G$ over $F$ by $\G^\un_{/F}$.

The unipotent completion of a finitely generated group is always defined over
$\Q$. Suppose that $\G$ is finitely generated.  Denote its unipotent completion
$\G^\un_{/\Q}$ over $\Q$ by $\G^\un$. Since $\big(\G^\un\otimes_\Q F\big)(F) =
\G^\un(F)$, the homomorphism $\G \to \G^\un(F)$ induces a homomorphism
$\G^\un_{/F} \to \G^\un\otimes_\Q F$. This homomorphism is an isomorphism.
Proofs can be found in \cite{hain:comp} and \cite{hain-matsumoto:p1}.

Denote the $n$th term of the lower central series of $\G^\un$ by $L^n \G^\un$.
The following ``well-known'' result follows from \cite{quillen} and
\cite{passman}, as explained in \cite[(2.17)]{hain-zucker}.

\begin{proposition}
\label{prop:lcs_unipt_compln}
If $\G$ is a finitely generated group and $F$ is a field of characteristic $0$,
then $D^n \G$ equals the inverse image of $(L^n\G^\un)(F)$ under the natural
homomorphism $\G \to \G^\un(F)$. Moreover, the image of $D^n\G$ is Zariski
dense in $L^n\G^\un$ and there is a natural isomorphism
$$
\big(D^n\G/D^{n+1}\big)\otimes_\Z F \cong (L^n \G^\un)(F)/(L^{n+1} \G^\un)(F).
$$
\end{proposition}

\subsection{Relative unipotent completion of mapping class groups}

Suppose that $S$ is a compact oriented surface of genus $g\ge 0$, that $r$ and
$n$ are non-negative integers such that $2g-2+n+r>0$, that $\x$ is a finite
subset of $S$ of cardinality $n$, and that $\vecs$ is a finite set of tangent
vectors of $S$ of cardinality $r$ as in Section~\ref{sec:mcgs}. The natural
homomorphism
$$
\rho : \G_{S,\x,\vecs} \to \Sp(H_1(S;F))
$$
is Zariski dense. Denote the relative unipotent completion of $\G_{S,\x,\vecs}$
with respect to $\rho$ by $\cG_{S,\x,\vecs}$. The relative completion of
$\G_{g,n,\vec{r}}$ will be denoted by $\cG_{g,n,\vec{r}}$.

Note that, when $g=0$ and $n+r\ge 3$, $H_1(S)$ is trivial so that relative
unipotent completion is unipotent completion: $\cG_{0,n,\vec{r}} =
\G_{0,n,\vec{r}/F}^\un$.
\end{example}

We recall some basic facts about relative unipotent completion of mapping class
groups. First, for all $g\ge 0$ and $n+r>0$, the natural homomorphism
$\G_{g,n,\vec{r}} \to \cG_{g,n,\vec{r}}$ is injective.

For a non-negative integer $m$, the level $m$ subgroup $\G_{g,n,\vec{r}}[m]$ of
$\G_{g,n,\vec{r}}$ is defined to be the kernel of the mod $m$ reduction $\rho_m
: \G_{g,n,\vec{r}} \to \Sp_g(\Z/m)$ of $\rho$. Note that
$\G_{g,n,\vec{r}}=\G_{g,n,\vec{r}}[1]$ and $T_{g,n,\vec{r}} =
\G_{g,n,\vec{r}}[0]$. When $g\ge 3$, the relative unipotent completion is
independent of the level.

\begin{theorem}[Hain \cite{hain:torelli}]
If $g\ge 3$, then for all non-negative integers $r,n$ and all positive integers
$m$, the inclusion $\G_{g,n,\vec{r}}[m] \to \G_{g,n,\vec{r}}$ induces an
isomorphism on relative unipotent completions. When $g=1,2$, the relative
completion of $\G_{g,n,\vec{r}}[m]$ depends non-trivially on the level $m>0$.
\qed
\end{theorem}

Right exactness implies that the sequence
$$
T_{g,n,\vec{r}}^\un \to \cG_{g,n,\vec{r}} \to \Sp_g \to 1
$$
is exact.

\begin{theorem}[Hain \cite{hain:comp,hain:torelli}]
\label{thm:cent_extn}
If $g\ge 3$, then the the kernel of $T_{g,n,\vec{r}}^\un \to \cG_{g,n,\vec{r}}$
is isomorphic to the additive group $\Ga$. That is, there is an exact sequence
$$
0 \to \Ga \to T_{g,n,\vec{r}}^\un \to \cG_{g,n,\vec{r}} \to \Sp_g \to 1. \qed
$$
\end{theorem}

\section{Unipotent Completion versus Pro-$\l$ Completion}

In this section we show that when $\G$ is a finitely generated $\l$-regular
group, the natural homomorphism $\G^\prol \to \G^\un(\Ql)$ is injective, a fact
we shall need in the proof of Theorem~\ref{thm:mcg_ladic}.

\begin{lemma}[{\cite[Lemma~A.7]{hain-matsumoto:p1}}]
If $\G$ is a finitely generated discrete group and $\ell$ a prime number, then
the natural homomorphism $\G \to \G^\un(\Ql)$ induces a homomorphism $\G^\prol
\to \G^\un(\Ql)$.
\end{lemma}

Next we give conditions under which this homomorphism is injective. In these
cases, the unipotent completion of a group $\G$ over $\Ql$ can be used to
``compute'' pro-$\l$ completion of $\G$.

As in previous sections, the lower central series of any group (discrete,
profinite, proalgebraic) $G$ will be denoted
$$
G = L^1G \supseteq L^2G \supseteq L^3G \supseteq \cdots
$$

\begin{proposition}
\label{prop:filt_D}
If $\G$ is a finitely generated $\l$-regular group, then the natural
homomorphism $\G^\prol \to \G^\un(\Ql)$ is injective.
\end{proposition}

\begin{proof}
The unipotent completion of any group is the inverse limit of the unipotent
completions of its torsion free nilpotent quotients
$$
\G^\un = \varprojlim_n (\G/D^n)^\un.
$$
By Proposition~\ref{prop:comp_rtfn}, the $\l$-regularity of $\G$ implies that
$$
\G^\prol \cong \varprojlim_{n} (\G/D^n)^\prol.
$$
It therefore suffices to prove the result when $\G$ is a finitely generated,
torsion free nilpotent group. This is easily proved by induction on the length
of the filtration $D^\dot$ of $\G$. If $D^2\G$ is trivial, then $\G$ is a
finitely generated abelian group, and thus isomorphic to $\Z^N$ for some $N$.
In this case
$$
\G^\prol \cong \G\otimes_\Z \Zl \text{ and }\G^\un(\Ql) \cong \G\otimes_\Z \Ql.
$$
So the result follows. If $D^{n+1}$ is trivial, then we can write $\G$ as an
extension
$$
1 \to D^n \G \to \G \to \G' \to 1
$$
where $\G' = \G/D^n$. Both $\G'$ and $D^n\G$ are finitely generated, torsion
free, nilpotent groups. By induction, we may assume that the result is true for
$\G'$ and $D^n\G$. To complete the proof, consider the commutative diagram
$$
\xymatrix{
1 \ar[r] & D^n\G \ar[r]\ar[d] & \G \ar[r]\ar[d] & \G' \ar[r]\ar[d] & 1 \cr
& (D^n\G)^\prol \ar[r]\ar[d]_{j_n} & \G^\prol
\ar[r]\ar[d]_j & (\G')^\prol \ar[r]\ar[d]_{j'} & 1 \cr
1 \ar[r] & (D^n\G)^\un(\Ql) \ar[r]^(.55)i & \G^\un(\Ql) \ar[r] & (\G')^\un(\Ql)
\ar[r] & 1 \cr
}
$$
The third row is exact by Proposition~\ref{prop:lcs_unipt_compln}; the second
is exact by the right exactness of pro-$\l$ completion. The maps $j_n$ and $j'$
are injective by induction. Since $i$ is injective, $j$ is as well.
\end{proof}

\section{Proof of Theorem~\ref{thm:mcg_ladic}}
\label{sec:proof}

Fix integers $g\ge 3$, $n\ge 0$, $r\ge 0$ and $m\ge 1$. Set $\G =
\G_{g,n,\vec{r}}$, $\G[m] = \G_{g,n,\vec{r}}[m]$ and $T = T_{g,n,\vec{r}}$.
Denote the $\Z$-group scheme $\Sp_g$ by $S$.

Fix a prime number $\ell$. Denote the pro-$\ell$ completion of $T$ by $T^\prol$
and the relative pro-$\ell$ completion of $\G[m]$ with respect to the standard
homomorphism to $S(\Zl)[m]$ by $\G[m]^\prol$. By Proposition~\ref{prop:S_ladic},
the relative pro-$\ell$ completion of $S(\Z)[m]$ is $S(\Zl)[m]$.

\subsection{An upper bound} In this section, we prove that if $T$ is
$\l$-regular, then the kernel of $T^\prol\to \G^\prol$ is either $0$ or
isomorphic to $\Zl$.

The right exactness of relative pro-$\l$ completion implies that the sequence
$$
T^\prol \to \G[m]^\prol \to S(\Zl)[m] \to 1
$$
is exact. Denote the kernel of $T^\prol  \to \G[m]^\prol$ by $K_{\l,m}$. When
$T$ is $\l$-regular,  we can give a tight upper bound on the size of $K_{\l,m}$.

\begin{lemma}
The kernel $K_{\l,m}$ of $T^\prol \to \G[m]^\prol$ is an abelian pro-$\l$ group
contained in the center of $T^\prol$. It has a natural $S(\Zl)[m]$ action.
\end{lemma}

\begin{proof}
By Proposition~\ref{prop:induced_action}, the conjugation action $\G \to \Aut T$
induces a continuous action $\G[m]^\prol \to \Aut T^\prol$. The kernel
$K_{\l,m}$ of $T^\prol \to \G[m]^\prol$ is thus a $\G[m]^\prol$-module that is
contained in the center of $T^\prol$. Since $T^\prol$ acts trivially on
$K_{\l,m}$, the $\G[m]$-action on $K_{\l,m}$ descends to a $S(\Zl)[m]$-action.
\end{proof}

We suspect, but cannot prove, that the $S(\Z)[m]$-action on $K_{\l,m}$ is
trivial. When $T$ is $\l$-regular, we can prove this and give a strong upper
bound on the size of $K_{\l,m}$.

\begin{proposition}
If $T$ is $\l$-regular, then the kernel $K_{\l,m}$ of $T^\prol \to \G[m]^\prol$
is a compact subgroup of $\Ql$ and is therefore either $0$ or isomorphic to
$\Zl$. The $S(\Z)[m]$-action on $K_{\l,m}$ is trivial.
\end{proposition}

\begin{proof}
Since $T$ is $\l$-regular, $j : T^\prol \to T^\un(\Ql)$ is injective. The
diagram
$$
\xymatrix{
&& T^\prol \ar[r]\ar[d]_j & \G[m]^\prol \ar[r]\ar[d] & S(\Zl)[m]\ar[d]\ar[r] & 1
\cr
0 \ar[r] & \Ql \ar[r] & T^\un(\Ql) \ar[r] & \cG(\Ql) \ar[r] & S(\Ql) \ar[r] & 1
}
$$
commutes and has exact rows by Theorem~\ref{thm:cent_extn} and the right
exactness of relative pro-$\l$ completion. Consequently, the kernel of $T^\prol
\to \G[m]^\prol$ is contained in the intersection of $\im j$ with $\Ql$, which
is compact as $T^\prol$ is profinite. The first assertion follows as the compact
subgroups of $\Ql$ are $0$ and $\ell^n\Zl$ where $n\in \Z$. The triviality of
the $S(\Z)[m]$-action on $K_{\l,m}$ follows from the fact that the inclusion
$\G[m] \to \G$ induces an isomorphism on relative unipotent completions
\cite[Prop.~3.3]{hain:torelli}, and because the kernel of $T^\un \to \cG$ is a
trivial $\cG$-module with respect to the action of $\cG$ on $T^\un$ that is
induced by the conjugation action of $\G$ on $T$.
\end{proof}

All remaining sections of Section~\ref{sec:proof} are devoted to proving that,
regardless of whether $T$ is $\l$-regular or not, the kernel of $T^\prol \to
\G^\prol$ contains copy of $\Zl$.

\subsection{The groups $G$ and $\hat{G}$.}

The argument to prove the non-triviality of the kernel of $T^\prol \to \G^\prol$
is somewhat intricate and uses the structure of some natural quotients of $\G$.

Consider the filtration
$$
\G = D^0\G \supseteq D^1 \G \supseteq D^2 \G \supseteq \cdots
$$
of
$\G$ defined, when $k\ge 1$, by $D^k \G := D^k T$. Since $T$ is finitely
generated \cite{johnson:fg}, the graded quotients $D^1/D^{k+1}$ are finitely
generated, torsion free nilpotent groups for all $k\ge 1$. Each graded quotient
is an $S(\Z)$-module.

Set
$$
V = D^1\G/D^2\G = H_1(T;\Z)/\text{torsion}.
$$
Then for all positive integers $m$ we have an extension
\begin{equation}
\label{eqn:extn}
1 \to V \to \G[m]/D^2 \to S(\Z)[m] \to 1
\end{equation}
Denote the relative pro-$\ell$ completion of $\G[m]/D^2$ with respect to the
natural homomorphism to $S(\Zl)[m]$ by $(\G[m]/D^2)^\prol$. Set $V_\ell =
V\otimes \Zl$.

\begin{lemma}
\label{lem:phi}
For all positive integers $m$, there is an injection
$$
\phi_m : \G[m]/D^2 \to S(\Z)[m]\ltimes V
$$
that commutes with the projections to $S(\Z)[m]$ and
whose restriction to $V$ is multiplication by $2$. Its image has finite
index.
\end{lemma}

\begin{proof}
It suffices to prove the case $m=1$. The extension (\ref{eqn:extn}) is
determined by a class $e\in H^2(S(\Z),V)$. Johnson's Theorem implies (cf.\
\cite{hain:torelli}) that $V$ is isomorphic to the sum of $\Lambda^3 H/H$ and
$n+r$ copies of $H$, where $H$ denotes the defining representation of $\Sp_g$.
Since $-I \in S(\Z)$ acts as $-1$ on $V$ and is central in $S(\Z)$, it follows
from the ``center kills'' argument that $H^2(S(\Z),V)$ is annihilated by $2$. In
particular, $2e=0$. Consequently, the extension of $S(\Z)$ by $V$ obtained by
pushing out the extension (\ref{eqn:extn}) along the multiplication by $2$ map
$\times 2 : V \to V$ splits. This implies that there is a map of extensions
\begin{equation}
\label{eq:diag}
\xymatrix{
1 \ar[r] & V \ar[r]\ar[d]^{\times 2} & \G/D^2 \ar[r]\ar[d]^{\phi} &
S(\Z)\ar@{=}[d]\ar[r] & 1
\cr
1 \ar[r] & V \ar[r] & S(\Z)\ltimes V \ar[r] & S(\Z) \ar[r] & 1.
}
\end{equation}
The result follows.
\end{proof}

\begin{corollary}
If $\ell$ is odd, the relative pro-$\ell$ completion of $\G[m]/D^2$ is
isomorphic to $S(\Zl)[m]\ltimes V_\ell$. When $\ell=2$,  $(\G[m]/D^2)^{(2)}$ has
finite index in $S(\Z_2)[m]\ltimes V_2$. \qed
\end{corollary}

\begin{proof}
We first show that the relative pro-$\l$-completion of $S(\Z)[m]\ltimes V$ is
$S(\Zl)[m]\ltimes V_\l$. The relative pro-$\l$ completion of $S(\Z)[m]$ is
$S(\Zl)[m]$ and the pro-$\l$ completion of $V$ is $V_\l$. The obvious
homomorphism
$$
\phi : S(\Z[m])\ltimes V \to S(\Zl)[m]\ltimes V_\ell
$$
is continuous. It therefore induces a homomorphism $\phi_\l$ such that the
diagram
$$
\xymatrix{
& V_\l \ar[r] \ar@{=}[d] & \big(S(\Z)\ltimes V\big)^\prol \ar[d]_{\phi_\l}
\ar[r] & S(\Zl)[m]\ar@{=}[d]\ar[r] & 1 \cr
1 \ar[r] & V_\l \ar[r] & S(\Zl)[m]\ltimes V_\ell \ar[r] & S(\Zl)[m] \ar[r] & 1
}
$$
commutes and has exact rows. It follows that $\phi_\l$ is an isomorphism.

When $\ell$ is odd, $2$ is a unit in $\Zl$. The result now follows from right
exactness of relative pro-$\l$ completion by taking the relative pro-$\l$
completion of the diagram
$$
\xymatrix{
1 \ar[r] & V \ar[r]\ar[d]^{\times 2} & \G[m]/D^2 \ar[r]\ar[d]^{\phi} &
S(\Z)[m]\ar@{=}[d]\ar[r] & 1
\cr
1 \ar[r] & V \ar[r] & S(\Z)[m]\ltimes V \ar[r] & S(\Z)[m] \ar[r] & 1.
}
$$
which is the restriction of the diagram (\ref{eq:diag}) to the level $m$
subgroups.
\end{proof}

The next step is to enlarge the quotient $\G[m]/D^2$ of $\G[m]$ by a central
extension. Since $V\otimes\Q$ is a rational representation of $S$ and since
$(\Gr^\dot_D T)\otimes\Q$ is a graded Lie algebra in the category of
$S(\Z)$-modules that is generated by $(\Gr^1_D T)\otimes \Q \cong V\otimes \Q$,
each $(D^kT/D^{k+1})\otimes\Q$ is a rational representation of the algebraic
group $S\otimes\Q$. By \cite[Thm.~10.1]{hain:torelli}, when $r=n=0$,
$$
\dim \big[\big(D^2T_g/D^3 T_g\big)\otimes \Q\big]^{S(\Z)} = 1.
$$
Let $\alpha_0 : D^2T_g/D^3 T_g \to \Q$ be an $S(\Z)$-invariant projection. Since
the action of $S(\Z)$ is semi-simple, this projection is unique up to scalar
multiplication. Define $\alpha_\Q : D^2T/D^3 \to \Q$ be the composite
$$
D^2T/D^3 \to D^2T_g/D^3 T_g \overset{\alpha_0}{\longrightarrow} \Q.
$$
Since $D^2T/D^3$ is finitely generated, the image of $\alpha_\Q$ is isomorphic
to $\Z$. Let $\alpha : D^2T/D^3 \to \Z$ be the unique (up to sign) surjection
whose composition with $\Z\to\Q$ is $\alpha_\Q$.

\subsubsection{Construction of $G$}
This is a key construction in the proof. Define $G$ to be the group obtained by
pushing out the extension
$$
0 \to D^2T/D^3 \to \G/D^3 \to \G/D^2 \to 1
$$
along the surjection $\alpha :  D^2T/D^3 \to \Z$. It is an extension of $\G/D^2$
by $\Z$. It can also be written as an extension
$$
1 \to E \to G \to S(\Z) \to 1,
$$
where the kernel $E$ is a central extension of $V$ by $\Z$, whose structure we
compute below.

\subsubsection{Heisenberg groups}
Suppose that $W$ and $A$ are abelian groups, and that $b : W\otimes_\Z W \to A$
is a bilinear pairing. The Heisenberg group $\Heis_b(W,A)$ is the group whose
underlying set is $W\times A$ and whose multiplication is
$$
(w,a)\cdot (w',a') = \big(w+w',a+a'+b(w,w')\big).
$$
This is a central extension
\begin{equation}
\label{eqn:heis}
0 \to A \to \Heis_b(W,A) \to W \to 0
\end{equation}
of $W$ by $A$. The commutator induces the skew symmetric bilinear pairing
$$
c : W\otimes_\Z W \to A
$$
defined by $c(w,w') = b(w,w') - b(w',w)$.\footnote{Note that if $b$ is skew
symmetric, then $c=2b$.} The abelianization of $\Heis_b(W,A)$ is thus an
extension
\begin{equation}
\label{eqn:abelianization}
0 \to A/\im c \to H_1(\Heis_b(W,A)) \to W \to 0.
\end{equation}

When $W$ is a free abelian group of finite rank, $c$ can be regarded as an
element of $H^2(W;A)$. It is the class of the extension (\ref{eqn:heis}) and
determines it up to isomorphism of extensions. Thus, when $W$ is torsion free of
finite rank, every central extension of $W$ by $A$ is of the form $\Heis_b(W,A)$
for some bilinear pairing $b : W\otimes W \to A$.

\begin{example}
The classical integral Heisenberg group is the group $\Heis_b(\Z^2,\Z)$, where
$b\big((x,y),(x',y')\big) = xy'$. It is isomorphic to the group
$$
\begin{pmatrix}
1 & \Z & \Z \cr 0 & 1 & \Z \cr 0 & 0 & 1
\end{pmatrix}
$$
The class of the extension is the class of the skew symmetrization
$\big((x,y),(x',y')\big) = xy'-x'y$ of $b$, which is the generator of
$H^2(\Z^2;\Z) \cong \Z$.
\end{example}

\begin{lemma}
\label{lem:heis_comp}
If $b: W\times W \to \Z$ is a bilinear pairing on a finitely generated, torsion
free abelian group, then the natural homomorphism
$$
\Heis_b(W,\Z) \to \Heis_b(W\otimes\Zl,\Zl)
$$
is the pro-$\l$ completion of $\Heis_b(W,\Z)$.
\end{lemma}

\begin{proof}
Consider the diagram
$$
\xymatrix{
& \Zl \ar[r]\ar@{=}[d] & \Heis_b(W,\Z)^\prol \ar[d]\ar[r]
& W\otimes \Zl \ar[r]\ar@{=}[d] & 1 \cr
1 \ar[r] & \Zl \ar[r] \ar[r] & \Heis_b(W\otimes\Zl,\Zl) \ar[r] &
W\otimes \Zl \ar[r] & 1
}
$$
Right exactness of pro-$\l$ completion implies that the first row is exact.
Since the first and third vertical maps are isomorphisms, the central map is also
an isomorphism. 
\end{proof}

Since $V$ is a finitely generated, torsion free group, the kernel $E$ of the
homomorphism $G \to S(\Z)$ is a Heisenberg group. The commutator of $E$ induces
a skew symmetric bilinear pairing $V\times V \to \Z$ which is $S(\Z)$-invariant
because it is invariant under the conjugation action of the mapping class group
on $T$.  The commutator pairing can be written as $dq$, where $d$ is a positive
integer and $q:V\times V \to \Z$ is a primitive pairing. (That is, it cannot be
divided by any integer $>1$.) In other words, the class of the extension is $dq
\in H^2(V,\Z)$. It follows from the sequence (\ref{eqn:abelianization}) that the
abelianization of $E$ is an extension
$$
0 \to \Z/d \to H_1(E) \to V \to 0.
$$

\begin{proposition}
\label{prop:eheis}
The integer $d$ divides $2$ and $E$ is a subgroup of $\Heis_q(V,\Z)$ of index
$2/d$.
\end{proposition}

\begin{proof}
Since $T \to E$ is surjective, $H_1(T) \to H_1(E)$ is surjective. The definition
of the rational dimension subgroups implies that $D^2T$ surjects onto the
torsion subgroup of $H_1(T)$. It follows that the torsion subgroup of $H_1(T)$
surjects onto $\Z/d$. Since the torsion subgroup of $H_1(T)$ has exponent $2$
\cite{johnson:h1}, it follows that $d=1$ or $2$. Standard facts about group
cohomology imply that if $d=2$, then $E\cong \Heis_q(V,\Z)$. If $d=1$, then $E$
has index 2 in $\Heis_q(V,\Z)$.
\end{proof}

Since $2$ is a unit in $\Zl$ when $\l$ is odd, Lemma~\ref{lem:heis_comp} and the
previous result give the following computation of $E^\prol$.

\begin{corollary}
\label{cor:Eprol}
If $\l$ is an odd prime number, then the pro-$\l$ completion of $E$ is
isomorphic to $\Heis_q(V\otimes\Zl,\Zl)$. If $\l=2$, then either $E^{(2)}$ is
isomorphic to $\Heis_q(V\otimes\Z_2,\Z_2)$ or has index 2 in it. In either case,
the completion $\Zl\to E^\prol$ of the inclusion $\Z\hookrightarrow E$ of the
central $\Z$ into $E$ is injective.
\end{corollary}

\subsubsection{Construction of $\hat{G}$} It will be necessary to compute the
relative pro-$\ell$ completion of $G$ up to a finite group. In order to do this,
we will construct another group $\hat{G}$ that contains $G$ as a finite index
subgroup.

\begin{proposition}
There is a group $\hat{G}$ that is an extension
\begin{equation}
\label{eqn:cent_ext}
0 \to \Z \to \hat{G} \to S(\Z)\ltimes V \to 1
\end{equation}
that contains $G$ as a finite index subgroup. More precisely, there is a
diagram of extensions
$$
\xymatrix{
0 \ar[r] & \Z\ar[r]\ar[d]^{\times 8/d} & G \ar[r]\ar[d] &
\G/D^2 \ar[r]\ar[d]^{\phi_1} & 1 \cr
0 \ar[r] & \Z\ar[r] & \hat{G} \ar[r] & S(\Z)\ltimes V \ar[r] & 1
}
$$
where $d\in \{1,2\}$ is the integer defined before Proposition~\ref{prop:eheis} and
$\phi_1$ is the homomorphism constructed in Lemma~\ref{lem:phi}. The extension
of $S(\Z)$ by $\Z$ obtained by restricting the extension (\ref{eqn:cent_ext}) to
$S(\Z)$ is $-(8g+4)$ times the universal central extension $\widetilde{S(\Z)}$
of $S(\Z)$.
\end{proposition}

It will follow from the proof of Corollary~\ref{cor:cent_extn} that $d$ is, in
fact, $2$.

\begin{proof}
The proof uses the geometry of moduli spaces of curves. We first review the
construction of the biextension bundle; a detailed exposition can be found in
\cite{hain-reed:arakelov}. In this proof, all moduli spaces are defined over
$\C$ and are viewed as analytic orbifolds. The moduli space $\A_g$ (viewed as an
orbifold) is a model of the classifying space $BS(\Z)$ of $S(\Z)$. The bundle of
intermediate jacobians $\J$ over $\A_g$ associated to the representation
$$
V_0 := \Lambda^3 H/H \cong T_g/D^2T_g
$$
of $S(\Z)$, where $H$ denotes the defining representation of $S$, is a model of
the classifying space $B\big(S(\Z)\ltimes V_0\big)$. A splitting of the
extension is given by the zero section of $\J(V_0) \to \A_g$. There is a natural
line bundle $\B_0 \to \J(V_0)$ whose associated $\C^\ast$-bundle $\B_0^\ast$ is
a model of $B\big(S(\Z)\ltimes \Heis_{q}(V_0,\Z)\big)$, where $q$ is the
$S$-invariant bilinear pairing $V\times V\to \Z$ defined before
Proposition~\ref{prop:eheis}. The restriction of $\B_0$ to the zero section
$\A_g$ is trivial. A detailed exposition of the construction of $\B_0$ is given
in \cite{hain-reed:arakelov}.

The period mapping $\M_g \to \A_g$ lifts to a holomorphic mapping $\nu_0 : \M_g
\to \J(V_0)$, which is the normal function that takes the moduli point of the
curve $C$ to the point in the primitive intermediate jacobian of $\Jac C$
that corresponds to the algebraic cycle $C-C^-$ in $\Jac C$:
$$
\xymatrix{
& \J(V_0)\ar[d] \cr
\M_g \ar[r]\ar[ur]^{\nu_0} & \A_g
}
$$
However, it does not lift to a section of $\B_0^\ast \to \A_g$. The obstruction
is the first Chern class of $\nu_0^\ast \B_0$, which is $(8g+4)\lambda \in
H^2(\M_g;\Z)\cong H^2(\A_g;\Z)$, where $\lambda$ is the first Chern class of the
line bundle $\L$ over $\A_g$ corresponding to the universal central extension of
$S(\Z)$, pulled back to $\M_g$. This Chern class computation follows from
results of Morita \cite[(5.8)]{morita} and is proved directly in
\cite{hain-reed:chern}.

The Chern class computation implies that the pullback of
$\B_0\otimes\L^{\otimes(-(8g+4))}$ to $\M_g$ along $\nu_0$ is trivial. This
implies that $\nu_0 : \M_g \to \J(V_0)$ lifts to a section $\nutilde$ of the
$\C^\ast$-bundle $(\B_0\otimes\L^{\otimes(-(8g+4))})^\ast$:
$$
\xymatrix{
&& \M_g\ar[dll]_{\nutilde}\ar[dl]_{\nu_0} \ar[d] \cr
\big(\B_0\otimes\L^{\otimes(-(8g+4))}\big)^\ast \ar[r] & \J(V_0) \ar[r] & \A_g
}
$$
At this stage it is useful to reinterpret this statement group theoretically.
Set $\hat{G}_0 = \pi_1((\B_0\otimes\L^{\otimes(-(8g+4))})^\ast,\ast)$. This
is an extension
$$
0 \to \Z \to \hat{G}_0 \to S(\Z)\ltimes V_0 \to 1
$$
that is also an extension of $S(\Z)$ by the Heisenberg group $\Heis_q(V_0,\Z)$.
Its restriction to $S(\Z)$ (fundamental group of the zero section of $\J(V_0)$)
is $\L^{\otimes(-(8g+4))}$, which has class $-(8g+4)\lambda \in H^2(S(\Z);\Z)$.
The lift $\nutilde$ of the normal function $\nu_0$ induces a homomorphism $\G_g
\to \hat{G}_0$ whose image contains the center $\Z$; the induced homomorphism
$T_g \to V_0$ is twice the Johnson homomorphism and has image $2V_0$.

This induces a homomorphism $\psi : E \to \Heis_q(V_0,\Z)$. Noting that the commutator
pairing of $\Heis_q(V_0,\Z)$ is $2q$, we see that $\psi$ must be multiplication
by $8/d$ on the central $\Z$'s:
$$
\xymatrix{
\Lambda^2 V_0 \ar[r]^{\times 4} \ar[d]_{dq} & \Lambda^2 V_0 \ar[d]_{2q} \cr
\Z \ar[r]^{\times 8/d} & \Z
}
$$
Here the first column is the commutator pairing of $E$ and the second is the
commutator pairing of $\Heis_q(V_0,\Z)$. This completes the proof of the result
when $r=n=0$ by taking $\hat{G} = \hat{G}_0$.

By replacing $\G_{g,n,\vec{r}}$ by $\G_{g,n+r}$, we may assume that $r=0$. When
$n=1$, $V=\Lambda^3 H$. There is a family of intermediate jacobians
$\J(\Lambda^3 H)$ over $\A_g$, which also fibers over $\J(V_0)$:
$$
\J(\Lambda^3 H) \to \J(V_0) \to \A_g.
$$
This has fundamental group $S(\Z)\ltimes \Lambda^3 H$ and the restriction of
$\nu_1$ to $T_{g,1}$ induces twice the Johnson homomorphism $T_{g,1} \to
\Lambda^3 H$. The normal function lifts to a normal function $\nu_1 : \M_{g,1}
\to \J(\Lambda^3 H)$ that induces the homomorphism $\G_{g,1} \to S(\Z)\ltimes
\Lambda^3 H$ of Lemma~\ref{lem:phi}. It is the normal function that takes the
moduli point of the pointed curve $(C,x)$ to the point in the intermediate
jacobian of $H_3(\Jac C)$ determined by the algebraic cycle $C_x - C_x^-$ in
$\Jac C$.

This generalizes to the case $n>1$ as follows: Denote the projection $\Lambda^3
H \to V_0$ by $p$. Then, by \cite[Cor.~3]{hain:farb},
$$
V = \{(u_1,\dots,u_n) \in \big(\Lambda^3 H\big)^n: p(u_1) = \dots = p(u_n)\}.
$$
Define $J(V)$ to be fibered product of $n$ copies of $\J(\Lambda^3 H)$ over
$\J(V_0)$. It is a model of the classifying space of $S(\Z)\ltimes V$. Each of
the $n$ points determines a normal function $\nu_j : \M_{g,n} \to \J(\Lambda^3
H)$. These all project to the normal function
$$
\M_{g,n}\to \M_g \overset{\nu_0}{\longrightarrow} \J(V_0).
$$
They therefore define a section $\nu_n : \M_{g,n} \to \J(V)$. This induces the
homomorphism $\G_{g,n} \to S(\Z)\ltimes V$ of Lemma~\ref{lem:phi}.

Denote the projection $\J(V) \to \J(V_0)$ by $\pi$. Then one has a diagram
$$
\xymatrix{
& \pi^\ast\big(\B_0\otimes\L^{\otimes(-(8g+4))}\big)^\ast \ar[r]\ar[d] &
\big(\B_0\otimes\L^{\otimes(-(8g+4))}\big)^\ast \ar[d] \cr
& \J(V) \ar[r]^\pi \ar[d] & \J(V_0) \ar[d] \cr
\M_{g,n} \ar[ruu]^{\nutilde_n} \ar[ru]^{\nu_n} \ar[r] & \A_g \ar@{=}[r] & \A_g
}
$$
The group $\hat{G}$ is defined to be
$$
\hat{G} = \pi_1(\pi^\ast\big(\B_0\otimes\L^{\otimes(-(8g+4))}\big)^\ast,\ast).
$$
The homomorphism $\G_{g,n} \to \hat{G}$ is induced by $\nutilde_n$. Its image
is $G$.
\end{proof}

Denote the inverse image of $S(\Z)[m]$ in $G$ by $G[m]$ and by $\hat{G}[m]$ the
inverse image of $S(\Z)[m]$ under the natural projection $\hat{G} \to S(\Z)$.
Then $G[m]$ has finite index in $\hat{G}[m]$.

\begin{corollary}
\label{cor:cent_extn}
For all positive integers $m$, there is a subgroup $\tilde{S}$ of $G[m]$ whose
image $S'$ under the quotient mapping $G[m] \to S(\Z)[m]$ has finite index and
which is, up to $4$-torsion, $-(2g+1)$ times the restriction of the universal
central extension of $S(\Z)$ to $S'$. Moreover, $d=2$ and $E$ is isomorphic to
$\Heis_q(V,\Z)$.
\end{corollary}

\begin{proof}
Denote by $\widehat{S(\Z)}$ the restriction of the extension
(\ref{eqn:cent_ext}) to $S(\Z)$. This is $-(8g+4)$ times the universal central
extension of $S(\Z)$. Let $\tilde{S}$ be the intersection of $G[m]$ with
$\widehat{S(\Z)}$ in $\hat{G}$. The statement about the Chern class of the
extension follows as $\tilde{S}$ contains the center of $G[m]$, which has index
$8/d$ in the center of $\widehat{S(\Z)}$. Its Chern class is thus, up to
$8/d$-torsion, $d/8$ times the restriction of the Chern class $-(8g+4)\lambda$
of the extension $\widehat{S(\Z)}$ of $S(\Z)$. Since $\lambda$ generates
$H^2(\Sp_g(\Z),\Z)$, this forces $d=2$. The remaining assertion follows from
Proposition~\ref{prop:eheis}.
\end{proof}

Recall that $G[m]$ is an extension
$$
1 \to E \to G[m] \to S(\Z)[m] \to 1.
$$

\begin{corollary}
\label{cor:kernel}
For all prime numbers $\ell$, the kernel of the homomorphism
$$
E^\prol \to G[m]^\prol
$$
is isomorphic to $\Zl$.
\end{corollary}

\begin{proof}
Apply relative pro-$\l$ completion to the diagram
$$
\xymatrix{
0 \ar[r] & \Z \ar[r]^\alpha\ar[d]_\beta & \tilde{S} \ar[r]\ar[d] &
S' \ar[r]\ar[d] & 1 \cr
1 \ar[r] & E \ar[r]^(.43)\gamma & G[m] \ar[r] & S(\Z)[m] \ar[r] & 1
}
$$
to obtain the commutative diagram
$$
\xymatrix{
\Zl \ar[r]^{\alpha^\prol}\ar[d]_{\beta^\prol} & \tilde{S}^\prol \ar[r]\ar[d] &
S'^\prol \ar[r]\ar[d] & 1 \cr
E^\prol \ar[r]^(.45){\gamma^\prol} & G[m]^\prol \ar[r] &
S(\Z)[m]^\prol \ar[r] & 1
}
$$
whose rows are exact. Corollary~\ref{cor:Eprol} implies that $\beta^\prol$ is
injective. Since $\tilde{S}$ has index $4$ in $\widehat{S(\Z)}$,
Proposition~\ref{prop:cent-extn} implies that the image of $\alpha^\prol$ is
finite, so that $\ker \alpha^\prol \cong \Zl$. The commutativity of the diagram
implies that $\beta^\prol(\ker\alpha^\prol)$ is contained in $\ker\gamma^\prol$.
Since $\ker\gamma^\prol$ is contained in $\ker\{E^\prol\to V^\prol\}$, which is
isomorphic to $\Zl$, this gives the result.
\end{proof}

\subsection{The kernel of $T^\prol \to \G^\prol$ contains $\Zl$}

We have the exact sequence
$$
1 \to K \to T^\prol \to \G[m]^\prol \to S(\Zl)[m] \to 1
$$
where $K:= K_{\l,m}$. Denote the kernel of the mapping $\G \to G$ by $C$. This
is a subgroup of $T$. Observe that $T/C$ is isomorphic to $\Heis_q(V,\Z)$ and
that $\G[m]/C$ is isomorphic to $G[m]$.

Denote the closure of the image of $C$ in $T^\prol$ (resp., $\G[m]^\prol$) by
$\Cbar_T$ (resp., $\Cbar_\G$). Since $C$ is normal in $\G[m]$,  $\Cbar_\G$ is
normal in $\G[m]^\prol$.  Right exactness of relative pro-$\l$ completion
implies that $G[m]^\prol \cong \G[m]^\prol/\Cbar_\G$ and $(T/C)^\prol \cong
T^\prol/\Cbar_T$. It also implies that the sequence
$$
1 \to K/(K\cap \Cbar_T) \to (T/C)^\prol \to G[m]^\prol \to S(\Zl)[m]
\to 1
$$
is exact.  Corollary~\ref{cor:kernel} implies that $K/(K\cap \Cbar_T) \cong
\Zl$. It follows that $K$ contains a copy of $\Zl$.

\subsection{Genus $2$}
\label{sec:genus2}
When $g=2$, the homomorphism $T^\prol \to \G^\prol[m]$ has a large kernel. For
simplicity, we restrict ourselves to the case $n=r=0$. Non-injectivity in the
case $r+n> 0$ follows by a similar argument. Details are left to the reader.

Set
$$
\Tbar := \im\{T_2^\prol \to \G_2^\prol[m]\}.
$$
Observe that its abelianization, $H_1(\Tbar)$ is a pro-$\l$ group with a
continuous $\Sp_2(\Zl)[m]$ action.

Mess' computations \cite{mess} imply that, as an $\Sp_2(\Z)$-module,
$$
H_1(T_2) \cong \Ind_{C_2 \ltimes \SL_2(\Z)^2}^{\Sp_2(\Z)} \Z.
$$
Here the group $C_2\ltimes \SL_2(\Z)^2$ is the stabilizer in $\Sp_2(\Z)$ of the
decomposition of the integral homology of an abelian surface which is the
product of two elliptic curves. The cyclic group $C_2$ interchanges the two
factors, and the two copies of $\SL_2(\Z)$ act on the homology of the two
elliptic factors.

It is easy to produce $\l$-adic quotients of $H_1(T_2)$ that cannot be
quotients of $H_1(\Tbar)$. For example, choose a prime number $p$ that does not
divide $m\l$. Let
$$
M = \Ind_{C_2 \ltimes \SL_2(\F_p)^2}^{\Sp_2(\F_p)} \Zl.
$$
View this as an $\Sp_2(\Z)[m]$-module via the projection $\Sp_2(\Z)[m] \to
\Sp_2(\F_p)$, which is surjective as $p$ does not divide $m$.

\begin{proposition}
\label{prop:ref}
The group $M$ is not a continuous quotient of $H_1(\Tbar)$. Consequently,
$H_1(T^\prol) \to H_1(\Tbar)$ is not injective. 
\end{proposition}

\begin{proof}
If there were a quotient mapping $\pi : H_1(\Tbar) \to M$, it would have to be
$\Sp_2(\Z)[m]$-equivariant. Since $\Sp_2(\Z)[p]$ acts trivially on $M$,
$\Sp_2(\Zl)[pm]$ would have to act trivially on $M$. But this is impossible as
$\Sp_2(\Zl)[pm] = \Sp_2(\Zl)[m]$ and since $\Sp_2(\Z)[m] \to \Sp_2(\F_p)$ is
surjective.
\end{proof}

An improvement of the previous result and the following corollary were suggested
by the referee.

\begin{corollary}
\label{cor:ref}
The kernel of $T_2^\prol \to \G_2[m]^\prol$ contains a copy of $\Zl$. \qed
\end{corollary}

\begin{proof}
Set $K = \ker\{T^\prol\to \Tbar\}$. The previous result implies that the
composite $K \to T \to M$ is non-trivial. Since $M$ is a free $\Zl$-module, the
image of $K \to M$ is a non-zero free $\Zl$-module, which implies that $K$
contains at least one copy of $\Zl$.
\end{proof}

\section{Proof of Theorem~\ref{thm:galois}}
\label{sec:unram}

First we prove the case $r=0$, from which the case $r>0$ follows by a group
theoretic argument presented at the end of this section.

Suppose that $2g-2+n>0$.  Denote the moduli stack of smooth projective
$n$-pointed genus $g$ curves by $\M_{g,n/\Z}$ and its Deligne-Mumford
compactification by $\Mbar_{g,n/\Z}$ \cite{knudsen,knudsen2}.

Several times in this section we need to normalize a scheme in another scheme
which is not irreducible. Specifically, suppose that $Y' \to Y$ is a finite flat
representable $1$-morphism of Deligne-Mumford stacks.  Suppose that $X$ is an
integral scheme and that $X' \to Y$ is a morphism from a nonempty Zariski open
subset of $X$ to $Y$. Define the normalization of $X$ with respect to
$$
X'\to Y \leftarrow Y'
$$
to be the disjoint union of the normalizations of $X$ in the function fields of
the irreducible components of $X'\times_Y Y'$.

\subsection{Tangential base points}
\label{sec:base_points}

The notion of tangential base point was introduced by Deligne in
\cite{deligne:tangential}. A $\Q$-rational tangential base point of a connected
scheme $X$ yields an exact functor from the category of finite \'etale coverings
of $X$ to the category of finite \'etale coverings of $\Spec\Q$, and hence
induces a homomorphism $\pi_1(\Spec\Q) \to \pi_1(X)$.  Since we deal only with
fundamental groups, we identify a tangential base point with the corresponding
exact functor. 

We consider the $\Q$-rational tangential base points constructed as follows.
Let 
$$
C \to \Spec\Z[[q_1,q_2,\ldots,q_{3g-3+n}]]
$$
be the universal deformation of a maximally degenerate stable curve of type
$(g,n)$ over $\Z$, where the parameter $q_i$ corresponds to a smoothing of the
$i$-th double point. Such deformations are constructed in
\cite[Rem.~2.3.10]{ihara-naka}. 

Let
$$
[C]: \Spec \Z[[q_1,q_2,\ldots,q_{3g-3+n}]] \to \Mbar_{g,n/\Z}
$$
be the classifying map. The choice of the set of parameters $q_i$ determines a
$\Q$-rational base point of $\M_{g,n/\Z}$ as follows.\footnote{If $q_i$ is
replaced by $-q_i$ for example,  the exact functor, and hence the base point,
changes.}

For each positive integer $m$, set
$$
B = \Q[[q_1,q_2,\ldots,q_{3g-3+n}]] \text{ and }
B^{1/m} = \Q[[q_1^{1/m},q_2^{1/m},\ldots,q_{3g-3+n}^{1/m}]].
$$
Then $[C]$ induces a map
$$
[C]_m : \Spec B^{1/m}[(q_1\cdots q_{3g-3+n})^{-1}] \to \M_{g,n/\Q}.
$$
Suppose that $M \to \M_{g,n/\Q}$ is a finite \'etale covering and $\Mtilde_m$ be
the normalization of $B^{1/m}$ with respect to
$$
\Spec B^{1/m}[(q_1\cdots q_{3g-3+n})^{-1}] \to \M_{g,n/\Q} \leftarrow M.
$$

By Abhyankar's Lemma \cite[Expos\'e XIII]{SGA1}, there is a positive integer $m$
such that $\Mtilde_m$ is \'etale over $B^{1/m}$. By specializing each
$q_i^{1/m}$ to $0$, we obtain a finite etale cover of $\Spec\Q$. This covering
is independent of the choice of $m$ and defines an exact functor from the
category of finite \'etale covers of $\M_{g,n/\Q}$ to the finite \'etale covers
of $\Spec\Q$. We shall abuse notation and denote this tangential base point of
$\M_{g,n}$ by $[C]$.

The tangential base point $[C]$ induces a section of the short exact sequence
\begin{equation}
\label{eq:short-standard}
1 \to \pi_1(\M_{g,n}\otimes\Qbar) \to \pi_1(\M_{g,n}\otimes \Q) \to G_\Q \to 1
\end{equation}
and thus a Galois action
$$
\rho_{[C]}: G_\Q \to \Aut \pi_1(\M_{g,n}\otimes\Qbar).
$$
The choice of an imbedding $\Qbar \hookrightarrow \C$ determines an isomorphism
$$
\pi_1(\M_{g,n}\otimes_\Z \Qbar) \cong \G_{g,n}^\wedge
$$
of the geometric fundamental group of $\M_{g,n}$ with the profinite completion
of the mapping class group $\G_{g,n}$. We will make this identification. By
Proposition~\ref{prop:closure}, the relative pro-$\ell$ completion of
$\G_{g,n}^\wedge$ is isomorphic to $\G_{g,n}^\prol$.

Our next task is to lift the base point $[C]$ of $\M_{g,n}$ to a base point
$[C']$ of $\M_{g,n+1}$. The special fiber $C_0$ ($q_j=0$, all $j$) of the curve
$C$ is a stable curve, each of whose components is isomorphic to $\P^1$; each
$\P^1$ has exactly three distinguished points (i.e., singular points or marked
points). Denote the $n$ marked points of $C$ by $x_1,\dots,x_n$.  Define $C_0'$
to be the curve obtained by replacing the first node ($q_1=0$) by a copy of
$\P^1$, where $0,\infty \in \P^1$ are nodes. Take $x_{n+1}$ to be the point $1$
on the added $\P^1$. Then $C'_0$ is an $n+1$ marked stable curve whose moduli
point in $[C_0']$ in $\Mbar_{g,n+1}$ maps to the moduli point $[C_0]$ in
$\Mbar_{g,n}$ by forgetting $(n+1)$-st point.

Extend $C'_0$ to an $(n+1)$-marked curve $C'$ over
$\Z[[q_1,\ldots,q_{3g-3+n+1}]]$ using the deformation theory construction in
\cite{ihara-naka}. Label the parameters so that $q_{3g-3+n+1}$ is the parameter
associated to the double point at the $\infty$ on the attached $\P^1$, $q_1$ the
parameter associated to that at $0$ on the attached $\P^1$, and $q_j$ is the
parameter corresponding to the $j$th node of $C_0$ when $1<j\le 3g-3+n$. 

The curve $C'$ corresponds to a morphism $[C']: \Spec
\Z[[q_1,\ldots,q_{3g-3+n+1}]] \to \Mbar_{g,n+1}$. If we forget the $(n+1)$st
marked point from $C'$, we have a family of semi-stable curves. Its
stabilization is isomorphic to  $C \times \Spec \Z[[q_1,\ldots,q_{3g-3+n+1}]]$.
Thus, the following diagram 
$$
\xymatrix{
\Spec \Z[[q_1,\ldots,q_{3g-3+n+1}]] \ar[d]_{} \ar[r]^(.7){[C']} &
\Mbar_{g,n+1} \ar[d]^{}\cr
\Spec \Z[[q_1,\ldots,q_{3g-3+n}]]   \ar[r]^(.7){[C]} & \Mbar_{g,n}.
}
$$
commutes, where the left vertical map is obtained by $q_1 \mapsto
q_1q_{3g-3+n+1}$ and $q_i \mapsto q_i$ when $1 < i \leq 3g-3+n$.

This and the above construction of the exact functors give a homomorphism 
$$
\pi_1(\M_{g,n+1},[C'])\to \pi_1(\M_{g,n},[C])
$$
that preserves the section from $G_\Q$. On the other hand, the image of $[C']$
is in the fiber isomorphic to $C$,  which gives a tangential base point of the
smooth fiber locus $C_{\Z((q_1,\dots,q_{3g-3+n}))}$.\footnote{For a ring $R$ we
denote $R[[t_1,\dots,t_r]][t_1^{-1},\dots,t_r^{-1}]$ by $R((t_1,\dots,t_r))$.}
Let $\Omega$ be an algebraically closed field that contains $B^{1/m}$ for all
positive integers $m$. Then we have a short exact sequence
\begin{equation}
\label{eq:exact-n}
1 \to \pi_1(C_\Omega, [C'])\to \pi_1(\M_{g,n+1/\Q},[C'])\to
\pi_1(\M_{g,n/\Q}, [C]) \to 1.
\end{equation}
The left group is isomorphic to the profinite completion of $\Pi_{g,n}$, the
fundamental group of an $n$-punctured, genus $g$ reference surface. This short
exact sequence gives the universal monodromy representation
\begin{equation}
\label{eq:univ-monod}
\pi_1(\M_{g,n/\Q},[C]) \to \Out \Pi_{g,n}^\prol.
\end{equation}
Restricting to the geometric fundamental group gives the geometric monodromy
representation
$$
\G_{g,n}^\wedge \cong \pi_1(\M_{g,n/\Qbar},[C])\to \Out \Pi_{g,n}^\prol.
$$
Further restricting to the action on the abelianization of
$\pi_1^\prol(C_\Omega)$,  we obtain a representation $\G_{g,n}^\wedge \to
\Sp_g(\Zl)$. The $G_\Q$-actions  (induced from (\ref{eq:exact-n}) and the
tangential section) on $\G_{g,n}^\wedge$ and $\Sp_g(\Zl)$ are
$G_\Q$-equivariant. Thus, the functoriality of relative pro-$\ell$ completion
implies that there is a Galois action
$$
\rho_{[C]}^\prol: G_\Q \to \Aut \G_{g,n}^\prol.
$$

\begin{theorem}[Mochizuki-Tamagawa]
\label{th:tamagawa-mochizuki}
For all prime numbers $\l$ and all $(g,n)$ satisfying $2g-2+n>0$, the
representation $\rho_{[C]}^{\prol}$ is  unramified outside $\l$. That is, it
factors through $\pi_1(\Spec \Z[1/\ell])$:
$$
G_\Q \to \pi_1(\Spec \Z[1/\ell]) \to \Aut \G_{g,n}^\prol.
$$
\end{theorem}

The case $g=0$ is proved by Ihara \cite{ihara:spherebraid}.

\subsection{\'Etale coverings of moduli spaces}

Suppose that $A$ is a ring contained in an algebraically closed field $\Omega$
of characteristic 0, and hence that $A$ is an integral domain. We have
$\G_{g,n}^\wedge \cong \pi_1(M_{g,n/\Omega}) \to \pi_1(\M_{g,n/A})$. We say that
a finite \'etale covering $M \to \M_{g,n/A}$ is {\em geometrically
relative-$\ell$}, if the corresponding $\G_{g,n}^\wedge$-action on fibers
factors through $\G_{g,n}^\prol$. This is independent of  the choice of $A \to
\Omega$.  The category of such coverings is a Galois category. By restricting to
the geometrically relative-$\ell$ coverings, we obtain the quotient
\begin{equation}
\label{eq:short-standard-rell}
1 \to \G_{g,n}^\prol \to \pi_1(\M_{g,n}\otimes \Q)' \to G_\Q \to 1.
\end{equation}
of (\ref{eq:short-standard}).

The existence of suitable finite coverings of $\M_{g,n}$ that have a
compactification that is smooth over $\Z[1/\ell]$ was established by de~Jong and
Pikaart when $n=0$ for all $\l$ in \cite{dejong-pikaart}, when $n>0$ and $\l$ is
odd  by Boggi and Pikaart in \cite{boggi-pikaart}, and when $n>0$ and $\l=2$ by
Pikaart in \cite{pikaart}. Their results needed here are summarized in the
following statement.

\begin{proposition}
\label{prop:normal-crossing}
For all prime numbers $\l$ and all $(g,n)$ satisfying $2g-2+n>0$,  there is a
geometrically relative-$\ell$ Galois covering $M \to \M_{g,n}$ defined over
$\Z[1/\ell]$ that satisfies:
\begin{enumerate}
\setcounter{enumi}{-1}

\item $M$ is an integral scheme;

\item the normalization $\overline{M}$ of $\Mbar_{g,n}$ in the function field
of $M$ is proper and smooth over $\Z[1/\ell]$;

\item the boundary  $\overline{M} \setminus M$ is a relative normal crossing
divisor over $\Z[1/\ell]$;

\item  the ramification index of the covering $\overline{M} \to \Mbar_{g,n}$
along any irreducible component of the boundary is a power of $\ell$;

\item  the natural map $\G_{g,n}^\prol \to \Sp_g(\Z/\ell)$ factors through the
quotient of $\G_{g,n}^\prol$ corresponding to the Galois covering $M\otimes
\Qbar \to \M_{g,n}\otimes \Qbar$ for odd $\ell$; for $\ell=2$, the same
statement with $\Sp_g(\Z/4)$.

\end{enumerate}
\end{proposition}

We shall explain how this statement follows from their results. For $G$ a finite
quotient of $\Pi_{g,n}$ by a characteristic  subgroup, Deligne and Mumford
\cite{DM} introduced the moduli stack $_G\M_{g,n}$ over $\Z[1/|G|]$ of smooth
projective curves of type $(g,n)$ with a Teichm\"uller structure of level $G$.
It is the normalization over $\Z[1/|G|]$ of  the finite Galois covering
$_G\M_{g,n/\Q} \to \M_{g,n/\Q}$ corresponding to the kernel of
$\pi_1(\M_{g,n/\Q}) \to \Out G$. The construction over $\Spec \Z[1/|G|]$ works
equally well when $G$ is a quotient of $\Pi_{g,n}$ by a finite index subgroup
stabilized by  $\pi_1(\M_{g,n/\Q})$.

Proposition~\ref{prop:normal-crossing} is proved by considering $M=
{_G\M}_{g,n/\Z[1/\ell]}$ for a suitable choice of $G$. Specifically take $G$ to
be:
\begin{enumerate}

\item the quotient group of $\Pi_{g,0}$ by the normal subgroup generated by the
fourth term of its lower central subgroup and all $\l$th powers when $\l$ is odd
and $n=0$;

\item the quotient group of $\Pi_{g,0}$ by the normal subgroup generated by the
fourth term of its lower central subgroup and all fourth powers when $\l=2$ and
$n=0$;

\item the quotient $\Pi_{g,n}/W^3 \Pi_{g,n}$, where $W^3$ denotes the third
term of the weight filtration defined in \cite{boggi-pikaart}), by all $\l$th
powers when $\l$ is odd and $n>0$;

\item the quotient $\Pi_{g,n}/W^4 \Pi_{g,n}$, where $W^4$ denotes the fourth
term of the weight filtration defined in \cite{boggi-pikaart}), by all fourth
powers when $\l=2$ and $n>0$.

\end{enumerate}
For details, see de~Jong and Pikaart \cite[Thm~3.1.1(iii) and
Prop~2.3.6]{dejong-pikaart} for the first two cases,  Boggi and Pikaart
\cite[Prop.~2.6(ii) and Prop~2.8]{boggi-pikaart} in the third case, and Pikaart
\cite[Thm~3.3.1 (2) and Thm~3.3.3 (7)]{pikaart} in the fourth case. The extra
condition for the case $\ell=2$ in (4) in 
Proposition~\ref{prop:normal-crossing} is to assure that $M$ is a scheme, not
just a stack.

Let $p\neq \ell$ be a prime number. For a subring $A$ of the field $\Omega$,
$\M_{g,n/A}$ denotes the base change to $A$. Denote by $\Z_p^\ur$ the maximal
\'etale cover of $\Z_p$, and by $\Q_p^\ur$ its fraction field.  There is a
natural morphism
$$
\M_{g,n/\Qbar_p} \overset{f}{\longrightarrow} \M_{g,n/\Z_p^\ur}. 
$$
Denote the Galois categories of geometrically relative-$\ell$ coverings of
$\M_{g,n/A}$ by $\cC(\M_{g,n/A})$. Note that the fundamental group of
$\cC(\M_{g,n/\Qbar_p})$ is isomorphic to $\G_{g,n}^\prol$.

\begin{proposition}
\label{prop:equivalence}
The above $f$ induces an equivalence of categories
$$
f^\ast : \cC(\M_{g,n/\Z_p^\ur}) \to \cC(\M_{g,n/\Qbar_p}).
$$
\end{proposition}

\begin{proof}
Take a geometrically relative $\ell$-covering $M$ as in
Proposition~\ref{prop:normal-crossing}. Choose a connected component of $M
\otimes \Z_p^\ur$. It is a connected object of $\cC(\M_{g,n/\Z_p^\ur})$, which
we denote by $M_{\Z_p^\ur}$. Since $M \to \M_{g,n}[1/\ell]$ is etale and $p\neq
\ell$, the base change $M_{\Qbar_p}$ of $M_{\Z_p^\ur}$ is a connected object of
$\cC(\M_{g,n/\Qbar_p})$. Since the boundary of $\overline{M}$ is a relatively
normal crossing divisor over $\Z[1/\l]$,  so is the boundary of the Zariski
closure of $M_{\Z_p^\ur}$ in  $\overline{M}_{\Z_p^\ur}$ over $\Z_p^\ur$. The
Base Change Theorem \cite[Expos\'e XII]{SGA1} implies that $f$ induces an
equivalence between the Galois categories of the finite \'etale $\ell$-coverings
of $M_{\Z_p^\ur}$  and $M_{\Qbar_p}$. Since $M$ is a geometrically
relative-$\ell$ Galois cover of $\M_{g,n}$ defined over $\Z[1/\ell]$ that
satisfies condition (4) of Proposition~\ref{prop:normal-crossing}, the category
$\cC(M_A)$ of $\ell$-coverings of $M_A$ is the subcategory of $\cC(\M_{g,n/A})$
consisting of all the objects (and the morphisms) over $M_A$  ($A=\Qbar_p$ or
$\Z_p^\ur$).

Choose a geometric point $x_{\Qbar_p}$ of $M_{\Qbar_p}$. Let $x_{\Z_p^\ur}$ be
its image in $M_{\Z_p^\ur}$. Let $x_A'$  denote their image in $\M_{g,n/A}$, for
$A=\Qbar_p$ or $\Z_p^\ur$.  Then  $\pi_1(\cC(M_A),x_A)$ is a finite index
subgroup of $\pi_1(\cC(\M_{g,n/A}),x_A')$, and the quotient set is identified
with the fiber $F(x_A')$ of $M_A$ over $x_A'$ pointed at $x_A$. The cardinality
of $F(x_A')$ is the covering degree of $M_A \to \M_{g,n/A}$, hence is the same
for $A=\Qbar_p$ and $\Z_p^\ur$. We have the following commutative diagram
$$
\xymatrix{
\pi_1(\cC(M_{\Qbar_p}),x_{\Qbar_p}) \ar[d] \ar@{^{(}->}[r]
& \pi_1(\cC(\M_{g,n/\Qbar_p}),x_{\Qbar_p}') \ar[d]\ar[r] & F(x_{\Qbar_p}')\ar[d]
\cr
\pi_1(\cC(M_{\Z_p^\ur}),x_{\Z_p^\ur})\ar@{^{(}->}[r]
& \pi_1(\cC(\M_{g,n/\Z_p^\ur}),x_{\Z_p^\ur}') \ar[r] & F(x_{\Z_p^\ur}'),
}
$$
where, in each row, the right-hand set is the quotient of the middle group by
the left-hand group. The left-hand vertical arrow has been proved to be an
isomorphism. Since the right-hand mapping is bijective, the middle vertical
arrow is an isomorphism. The result follows. 
\end{proof}

\subsection{Proof of Theorem~\ref{th:tamagawa-mochizuki}}

It suffices to show that for every $p \neq \ell$, the inertia group
$\pi_1(\Q_p^\ur)$ at $p$ acts trivially on the relative pro-$\ell$ completion of
the mapping class group.

We use the notation of Section~\ref{sec:base_points}. Take $M$ over $\Z[1/\ell]$
as in Proposition~\ref{prop:normal-crossing} and $M_{\Z_p^\ur}$ as in the proof
of Proposition~\ref{prop:equivalence}. A tangential base point of $M$ over $[C]$
induces an exact functor from  $\cC(M_{\Z_p^\ur})$ to the category
$\cC(\Z_p^\ur)$ of finite \'etale coverings of $\Z_p^\ur$ as follows. 

Let $\Mtilde$ be the normalization of $\Spec\Z_p^{\ur}[[q_1,\ldots,q_{3g-3+n}]]$
with respect to
$$
\Spec\Z_p^{\ur}[[q_1,\ldots,q_{3g-3+n}]][(q_1\dots q_{3g-3+n})^{-1}]
\to \M_{g,n/\Z_p^\ur} \leftarrow M_{\Z_p^\ur}.
$$
The second and third assertions of Proposition~\ref{prop:normal-crossing}
concerning ramification imply that there are $\l$-power integers
$k_1,\ldots,k_{3g-3+n}$ (possibly 1) such that
$$
\Mtilde\cong\Spec\Z_p^{\ur}[[q_1^{1/k_1},\ldots,q_{3g-3+n}^{1/k_{3g-3+n}}]].
$$
This gives a $\Z_p^{\ur}$-rational tangential base point on $M_{\Z_p^{\ur}}$.
Abhyankar's Lemma \cite[Expos\'e XIII]{SGA1} implies that, for each finite
\'etale $\ell$-covering $N$ of $M_{\Z_p^\ur}$, there is an $\ell$-power integer
$j$ such that
$$
\widetilde{N} \to
\Spec\Z_p^{\ur}[[q_1^{1/jk_1},\ldots,q_{3g-3+n}^{1/jk_{3g-3+n}}]]
$$
is a finite \'etale morphism, where $\widetilde{N}$ is the normalization of
$$
\Spec\Z_p^{\ur}[[q_1^{1/jk_1},\ldots,q_{3g-3+n}^{1/jk_{3g-3+n}}]]
$$
with respect to
$$
\Spec\Z_p^{\ur}[[q_1^{1/jk_1},\ldots,q_{3g-3+n}^{1/jk_{3g-3+n}}]]
[(q_1\dots q_{3g-3+n})^{-1}]
\to M_{\Z_p^\ur} \leftarrow N
$$

By specializing the $q_i^{1/jk_i}$ to zero, we obtain a finite \'etale covering
of $\Z_p^\ur$. This gives an exact functor $\cC(M_{\Z_p^\ur}) \to
\cC(\Z_p^\ur)$. The geometric point $\xi: \Spec\Qbar_p \to \Spec \Z_p^\ur$ gives
a fiber functor $F_\xi$ of $\cC(\Z_p^\ur)$, whose composition with the above
functor is the fiber functor $F_\xi'$ of $\cC(M_{\Z_p^\ur})$ associated to 
$q_1^{1/k_1},\ldots,q_{3g-3+n}^{1/k_{3g-3+n}}$. Composing with the functor
$\cC(\M_{g,n/\Z_p^\ur})\to \cC(M_{\Z_p^\ur})$ gives a fiber functor $F_\xi''$ of
$\cC(\M_{g,n/\Z_p^\ur})$:
$$
\xymatrix{
\cC(\M_{g,n/\Z_p^\ur}) \ar[r]\ar[dr]_{F_\xi''} &
\cC(M_{\Z_p^\ur}) \ar[r]\ar[d]_{F_\xi'} &
\cC(\Z_p^\ur)\ar[dl]^{F_\xi}\cr
& {\mathsf{Sets}}
}
$$
Note that the fiber functor $[C]_*(\xi)$ equals $F_\xi''$. Define
$[C]_{\Z_p^\ur}$ to be the composite
$$
\pi_1(\Z_p^\ur, F_\xi) \to \pi_1(\cC(M_{\Z_p^\ur}), F_\xi')
\to \pi_1(\cC(\M_{g,n/\Z_p^\ur}),F_\xi'').
$$
This and the similar construction over $\Q_p^\ur$ gives a commutative diagram
$$
\xymatrix{
& \pi_1(\M_{g,n/\Qbar_p},F_\xi'')^\prol \ar@{^{(}->}[d]\cr
\pi_1(\Q_p^\ur, F_\xi)\ar[d] \ar[r]^(.4){[C]_{\Q_p^\ur}} \ar[r] &
\pi_1(\cC(\M_{g,n/\Q_p^\ur}),F_\xi'') \ar[d] \cr
\pi_1(\Z_p^\ur, F_\xi) \ar[r]^(.4){[C]_{\Z_p^\ur}} \ar[r] & 
\pi_1(\cC(\M_{g,n/\Z_p^\ur}),F_\xi'').
}
$$
where the top right-hand arrow is the inclusion of a normal subgroup. The top
right-hand group is isomorphic to $\G_{g,n}^\prol$. It is a normal subgroup of
$\pi_1(\cC(\M_{g,n/\Q_p^\ur}),F_\xi'')$. Proposition~\ref{prop:equivalence}
implies that the composite of the two right-hand vertical maps is an
isomorphism. The action of the inertia $\pi_1(\Q_p^\ur)$ at $p$ on
$\G_{g,n}^\prol$ is conjugation by an element of 
$\pi_1(\cC(\M_{g,n/\Q_p^\ur}),F_\xi'')$. The commutativity of the diagram
implies that the inertia action factors through $\pi_1(\Z_p^\ur, F_\xi)$, which
is trivial.  Theorem~\ref{th:tamagawa-mochizuki} follows.

Denote the kernel of $G_\Q \to \pi_1(\Z[1/\ell])$ by $J_\ell$. Define
$\G_{g,n}^{\arith, \prol}$ to be the quotient of $\pi_1(\M_{g,n/\Q},[C])$ by the
normal subgroup generated by $[C]_\ast(J_\ell)$ and the kernel of
$\Ghat_{g,n}\to \G_{g,n}^\prol$. Taking the quotient of the middle and the
right-hand groups in the short exact sequence (\ref{eq:short-standard-rell}), we
obtain the following result.

\begin{corollary}
\label{cor:exact_over_integer}
For all prime numbers $\l$ and all $(g,n)$ satisfying $2g-2+n>0$, there is a
split exact sequence
$$
1 \to \G_{g,n}^\prol \to \G_{g,n}^{\arith, \prol} \to \pi_1(\Z[1/\ell]) \to 1. 
$$
\end{corollary}

\begin{proof}
Consider the section of the short exact sequence (\ref{eq:short-standard-rell})
given by the tangential base point $[C]$. Theorem~\ref{th:tamagawa-mochizuki}
implies that the image of the kernel of $G_\Q \to \pi_1(\Z[1/\ell])$ under the
section induced by $[C]$ centralizes $\G_{g,n}^\prol$, hence is a normal
subgroup of the middle group. We can divide the middle and the right-hand groups
by this kernel to obtain the short exact sequence. The splitting is induced by
$[C]$.
\end{proof}

The following corollary was suggested by the referee.

\begin{corollary}
The group $\G_{g,n}^{\arith, \prol}$ is canonically isomorphic to the
fundamental group of the Galois category $\cC(\M_{g,n/\Z[1/\ell]})$ with respect
to the tangential base point $[C]$.  In particular, although the definition of
$\G_{g,n}^{\arith, \prol}$ depends on the choice of $[C]$, two different choices
yield two isomorphic groups. The isomorphism is canonical up to an inner
automorphism.  
\end{corollary}

\begin{proof}
Denote the kernel of $G_\Q \to \pi_1(\Z[1/\ell])$ by $J_\ell$. By the definition
of $\G_{g,n}^{\arith, \prol}$ given above, we have an identification
$$
\G_{g,n}^{\arith, \prol}= \pi_1\big(\cC(\M_{g,n/\Q})\big)/([C]_*(J_\ell)).
$$
On the other hand, by the ``purity of the branched locus'', we have
$$
\pi_1\big(\cC(\M_{g,n/\Z[1/\ell]})\big)
\cong \pi_1\big(\cC(\M_{g,n/\Q})\big)/\langle T_p|p\neq \ell\rangle,
$$
where $T_p$ denotes the inertia subgroup of  $\pi_1(\cC(\M_{g,n/\Q}))$ along the
divisor $\M_{g,n/\F_p}$ of $\M_{g,n/\Z[1/\ell]}$, and $\langle T_p|p\neq
\ell\rangle$ denotes the normal subgroup topologically generated by the union of
the $T_p$ for all primes $p \neq \ell$.  The result will follow from the
statement that these two quotients are identical as it implies that 
$\G_{g,n}^{\arith, \prol}$ is canonically  isomorphic to the fundamental group
of the Galois category $\cC(\M_{g,n/\Z[1/\ell]})$ with fiber functor $[C]$,  and
another choice of $[C]$ yields an isomorphic fundamental group, with isomorphism
determined up to an inner automorphism. The rest of the proof consists in
proving this equality.

When a covering in $\cC(\M_{g,n/\Z[1/\ell]})$ is pulled back along the
tangential section $[C]$ over $\Z[1/\ell]$,  it gives a finite \'etale covering
of $\Z[1/\ell]$. This implies  that $[C]_*(J_\ell)$ acts trivially on the
corresponding fiber,  and hence that $[C]_*(J_\ell) \subseteq \langle T_p|p\neq
\ell\rangle$. To establish the reverse inclusion, it suffices to show that the
image of $T_p$ in $\G_{g,n}^{\arith,\prol}$ is trivial.  Consider the exact
sequence
$$
1 \to \langle\tilde{T}_p\rangle \to \pi_1\big(\cC(\M_{g,n/\Q_p^\ur})\big) \to
\pi_1\big(\cC(\M_{g,n/\Z_p^\ur})\big) \to 1,
$$
where $\tilde{T}_p$ is the inertia subgroup of $\pi_1(\cC(\M_{g,n/\Q_p^\ur}))$,
and $\langle\tilde{T}_p \rangle$ is the normal subgroup topologically generated
by $\tilde{T}_p$.
The construction in Corollary~\ref{cor:exact_over_integer}, when $\Z[1/\ell]$ is
replaced with $\Z_p^\ur$, gives a  short exact sequence
$$
1 \to \G_{g,n}^\prol \to \pi_1\big(\cC(\M_{g,n/\Q_p^\ur})\big)/
[C]_*(\pi_1(\Q_p^\ur)) \to \pi_1(\Z_p^\ur) \to 1.
$$
Since $\pi_1(\Z_p^\ur)$ is trivial, the left group is isomorphic to the middle
group. Consider the commutative diagram
$$
\xymatrix{
\tilde{T}_p \ar@{->>}[d]_{} \ar[r] &
\pi_1\big(\cC(\M_{g,n/\Q_p^\ur})\big)/[C]_\ast(\pi_1(\Q_p^\ur)) \ar[d]\cr
T_p \ar[r] & \pi_1\big(\cC(\M_{g,n/\Q})\big)/[C]_\ast(J_\ell) \ar@{=}[r] &
\G_{g,n}^{\arith,\prol}.
}
$$
Since the left-hand vertical map is surjective, the triviality of the bottom
horizontal map will follow from the triviality of the upper horizontal map.
Since
$$
\xymatrix{
\pi_1(\Q_p^\ur) \ar[d]_{[C]_*} \ar[r] & \pi_1(\Z_p^\ur)
 \ar[d]^{[C]_*}\cr
\pi_1\big(\cC(\M_{g,n/\Q_p^\ur})\big) \ar[r] &
\pi_1\big(\cC(\M_{g,n/\Z_p^\ur})\big)
}
$$
commutes and the top right group is trivial, the bottom  horizontal map
factors through a homomorphism
$$
\pi_1\big(\cC(\M_{g,n/\Q_p^\ur})\big)/[C]_*(\pi_1(\Q_p^\ur)) 
\to  \pi_1\big(\cC(\M_{g,n/\Z_p^\ur})\big),
$$
which is an isomorphism because the both groups are isomorphic to
$\G_{g,n}^\prol$ by the above argument and Proposition~\ref{prop:equivalence}.
Since the image of $\tilde{T}_p$ is trivial in the right-hand group, it is also
trivial in the left-group.  This establishes the isomorphism of
$\G_{g,n}^{\arith, \prol}$ with the fundamental group of the Galois category
$\cC(\M_{g,n/\Z[1/\ell]})$ and with it, the result.
\end{proof}

\begin{theorem}
For all prime numbers $\l$ and all $(g,n)$ satisfying $2g-2+n>0$, the universal
monodromy representation (\ref{eq:univ-monod})
$$
\pi_1(\M_{g,n/\Q},[C]) \to \Out \pi_1^\prol(C_\Omega)
$$
factors through $\G_{g,n}^{\arith, \prol}$.
\end{theorem}

\begin{proof}
Since the universal monodromy representation is the outer representation
associated to the extension (\ref{eq:exact-n}), it suffices to show that
$$
1 \to \pi_1(C_\Omega)^\prol \to \G_{g,n+1}^{\arith, \prol} \to
\G_{g,n}^{\arith, \prol} \to 1
$$
is exact. This is equivalent to the exactness of
$$
1 \to \pi_1(C_\Omega)^\prol \to \G_{g,n+1}^{\prol} \to
\G_{g,n}^{\prol} \to 1,
$$
which was proved in Proposition~\ref{prop:injectivity}(2).
\end{proof}

\begin{remark} The proof of Theorem~\ref{th:tamagawa-mochizuki} was communicated
to us by Tamagawa. Mochizuki gave a more sophisticated proof, sketched below,
which uses log geometry and avoids the constructions of Boggi, de~Jong, and
Pikaart in Proposition~\ref{prop:normal-crossing}.

In the above proof, the problem of proving the result is reduced to  studying
the $\l$-coverings of $M$ in Proposition~\ref{prop:normal-crossing}, where
$M\subset \overline{M}$ is the complement of a normal  crossing divisor. This
allows the use of \cite[Expos\'e XII]{SGA1}. Mochizuki, on the other hand, 
takes $N \to \M_{g,n/\Z[1/\l]}$ to be the covering corresponding to
$\pi_1(\M_{g,n/\Z[1/\l]}) \to \GSp_g(\Z/\l)$.  He equips $\Mbar_{g,n}$ with the
standard log structure. The normalization $\overline{N}$ of
$\Mbar_{g,n/\Z[1/\l]}$ is then log-regular since the ramification is tame.
Mochizuki uses the log-purity theorem of Kato and Fujiwara (unpublished; see
\cite{mochizuki} for a proof) in place of \cite[Expos\'e XII]{SGA1} to complete
the proof.

\end{remark}

\subsection{The case $r>0$.}

Theorem~\ref{th:tamagawa-mochizuki} establishes Theorem~\ref{thm:galois} in the
case $r=0$. The case $r>0$ reduces to the case $r=0$.

We first discuss the Deligne-Mumford compactification of $\M_{g,n,\vec{r}}$.
Suppose that $2g-2+n+r>0$. Define a stable curve of type $(g,n,\vec{r})$ to be 
a stable curve of type $(g,n+r)$ together with the choice of a non-zero tangent
vector at each of the last $r$ points. The moduli space of stable curves of type
$(g,n,\vec{r})$ is a $(\Gm)^r$-bundle over $\Mbar_{g,n+r}$, which is not
complete if $r>0$. Denote the line bundle associated to the $j$th factor of this
$(\Gm)^r$ bundle by $L_j$. Its fiber over the stable curve
$[C,x_1,\dots,x_n,y_1,\dots,y_r]$ is $T_{y_j}C$, the tangent space of $C$ at
$y_j$. The $j$th factor of the $(\Gm)^r$ bundle over $\Mbar_{g,n+r}$ has fiber
$T_{y_j} C - \{0\}$ over $[C,x_1,\dots,x_n,y_1,\dots,y_r]$ and is therefore
defined over $\Spec \Z$. We compactify the $j$th factor of the $(\Gm)^r$ bundle
in the usual way; namely:
$$
\P\big(L_j \oplus \O_{\Mbar_{g,n+r}}\big).
$$
Define $\Mbar_{g,n,\vec{r}}$ to be the fibered product over $\Mbar_{g,n+r}$ of
the bundles $\P(L_j \oplus \O_{\Mbar_{g,n+r}})$. This is defined over $\Spec\Z$
as each of its factors is.

The construction of the deformation of maximally degenerate curves over $\Z$,
and the construction of $\Q$-rational tangential base points are all similar to
constructions in the case $r=0$ explained above. We start with a maximally
degenerate stable curve of type $(g,n+r)$ and then construct the universal
deformation over $\Z$,
$$
[C]: \Spec\Z[[q_1,q_2,\ldots,q_{3g-3+n+r}]] \to \Mbar_{g,n+r/\Z}. 
$$
The fiber product  $\Spec\Z[[q_1,q_2,\ldots,q_{3g-3+n+r}]] 
\times_{\Mbar_{g,n+r/\Z}} \Mbar_{g,n,\vec{r}}$ is $(\P^1)^r$-bundle over
$\Spec\Z[[q_1,q_2,\ldots,q_{3g-3+n+r}]]$. Take $q_1',\ldots,q_r'$ to be the
standard coordinate of  each copy of $\P^1$ so that $q_i'=0$ gives the zero of
the $i$-th copy of $\P^1\supset \Gm$. Now
$$ 
[C']: \Spec\Z[[q_1,q_2,\ldots,q_{3g-3+n+r}, q_1',\ldots,q_r']]
\to \Mbar_{g,n,\vec{r}/\Z}
$$
gives a tangential base point compatible with $[C]$. In the following, we use 
the Galois action given by these tangential base points.

Suppose that $r>0$. We have to establish the triviality of the action on
$\G_{g,n,\vec{r}}^\prol$ of the inertia subgroup $I_p$ of $G_\Q$ at $p$. By
Proposition~\ref{prop:injective}, the sequence
$$
0 \to \Zl(1)^r \to \G_{g,n,\vec{r}}^\prol \to \G_{g,n+r}^\prol \to 1
$$
is exact. Fix $\sigma \in I_p$ and take $g \in \G_{g,n,\vec{r}}^\prol$.  Then, 
since $I_p$ acts trivially on $\G_{g,n+r}^\prol$ (this is the case $r=0$), we
can define a function
$$
\varphi : \G_{g,n,\vec{r}}^\prol \to \Zl(1)^r
$$
by $\varphi: h \mapsto \sigma(h)h^{-1}$.  Since $\Zl(1)^r$ is central in
$\G_{g,n,\vec{r}}^\prol$,  $\varphi$ is a group homomorphism.  Since $I_p$ acts
trivially on $\Zl(1)^r$, $\varphi$ induces a homomorphism
$$
\phibar : \G_{g,n+r}^\prol \to \Zl(1)^r.
$$
Since $\G_{g,n+r}$ is finitely generated, $H^1(\G_{g,n+r},\Z)$ is a finitely
generated torsion free abelian group. Since $H^1(\G_{g,n+r},\Q)$ vanishes for
all $g\ge 1$,\footnote{This follows from Harer's work \cite{harer}. He proved
that $H_1(\G_g)$ vanishes when $g\ge 3$ and is cyclic of order 10 when $g=2$. It
is classical that $H_1(\G_{1,1})=\Z/12$. The result when $n+r>0$ follows by
induction and a standard and straightforward spectral sequence argument. It is
also proved in \cite[Prop.~5.1]{hain:normal} for $g\ge 3$ by different means.}
$H^1(\G_{g,n+r},\Z)$ vanishes whenever $g\ge 1$. This implies that $\phibar$ is
trivial whenever $g\ge 1$ as $\G_{g,n+r}$ is dense in $\G_{g,n+r}^\prol$. It
follows that $I_p$ acts trivially on $\G_{g,n,\vec{r}}^\prol$. When $g=0$,
relative pro-$\l$ completion coincides with pro-$\l$ completion, so that we can
appeal to \cite[Expos\'e XII]{SGA1}.

\end{document}